\documentclass[12pt]{article}

\usepackage{amsmath}
\usepackage{amssymb}
\usepackage{amsfonts}
\usepackage{amsthm} 
\usepackage{enumerate}

\newtheorem{theorem}{Theorem}[section]
\newtheorem{proposition}[theorem]{Proposition}
\newtheorem{lemma}[theorem]{Lemma}
\newtheorem{definition}[theorem]{Definition}

\newcommand{\rr}{{\mathbb R}}
\newcommand{\zz}{{\mathbb Z}}
\newcommand{\nn}{{\mathbb N}}

\newcommand{\ee}{{\mathbb E}}

\newcommand{\pp}{{\mathbb P}}

\newcommand{\calb}{{\cal B}}
\newcommand{\calee}{{\cal E}}
\newcommand{\calx}{{\cal X}}

\newcommand{\dist}{\hbox{dist\,}}
\newcommand{\supp}{\hbox{supp\,}}

\newcommand{\psieps}{\psi_\epsilon}
\newcommand{\phieps}{\phi_\epsilon}
\newcommand{\whm}{\widehat{\mu}}

\begin{document}

\title{Arithmetic progressions in sets of fractional dimension}
\author{Izabella {\L}aba and Malabika Pramanik}
\date{January 10, 2008}
\maketitle

\begin{abstract}
Let $E\subset\rr$ be a closed set of Hausdorff dimension $\alpha$.
We prove that if $\alpha$ is sufficiently close to 1, and if $E$ supports
a probabilistic measure obeying appropriate dimensionality and
Fourier decay conditions, then $E$ contains non-trivial 3-term
arithmetic progressions.

\end{abstract}

\section{Introduction}

\begin{definition}\label{def-propE}
Let $A\subset\rr$ be a set.  We will say that $A$ is {\em universal}\footnote{We are using the
terminology of \cite{mk}.} for a class $\calee$
of subsets of $\rr$ if any set $E\in\calee$ must contain an
affine (i.e. translated and rescaled) copy of $A$.
\end{definition}

If $\calee$ is the class of all subsets of $\rr$ of positive Lebesgue measure, then
it follows from Lebesgue's theorem on density
points that every finite set $A$ is universal for $\calee$.  Namely, let $E$ have positive 
Lebesgue measure, then $E$ has density 1 at almost every 
$x\in E$.  In particular, given any $\delta>0$, we may choose an interval $I=(x-\epsilon,
x+\epsilon)$ such that $|E\cap I|\geq (1-\delta)|I|$.  If $\delta$ was chosen small
enough depending on $A$, the set $E\cap I$ will contain an affine copy of $A$.

An old question due to Erd\H{o}s \cite{erdos} is whether 
any {\em infinite} set $A\subset\rr$ can be universal for all sets of positive 
Lebesgue measure.
It is known that not all infinite sets are universal: for instance,
if $A=\{a_n\}_{n=1}^\infty$ is a slowly decaying sequence such that $a_n\to 0$ 
and $\frac{a_{n-1}}{a_n}\to 1$, then one can construct explicit Cantor-type sets of
positive Lebesgue measure
which do not contain an affine copy of $A$ \cite{falconer2}.  Other classes of counterexamples
are given in \cite{bourg-sequences}, \cite{mk}; see also \cite{HL}, \cite{komjath} for
further related work.  There are no known examples of infinite sets
$A$ which are universal for the class of sets of positive measure.  In particular, the 
question remains open for $A=\{2^{-n}\}_{n=1}^\infty$.

The purpose of this paper is to address a related question,
suggested to the first author by Alex Iosevich:
if $A\subset \rr$ is a finite set and $E\subset[0,1]$ is a set of Hausdorff 
dimension $\alpha\in[0,1]$, must $E$ contain an affine copy of $A$?  
In other words, are finite sets universal for the class of all sets of Hausdorff dimension 
$\alpha$?  This more general statement already fails
if $A=\{0,1,2\}$ and $E$ is a set of Hausdorff dimension 1 but Lebesgue measure 0.
This is due to Keleti \cite{keleti}, who actually proved a stronger result: there is 
a closed set $E\subset [0,1]$ of Hausdorff dimension $1$ such that $E$ does not contain
any ``rectangle" $\{x,x+r,y,y+r\}$ with $x\neq y$ and $r\neq 0$.

In light of the above, one may ask if there is a natural subclass of sets of fractional
dimension for which a finite set such as $\{0,1,2\}$ might be universal.
This question is addressed in Theorem \ref{thm-main}, which is the main result of this article.
We define the Fourier coefficients of a measure $\mu$ supported on $[0,1]$ as
$$
\widehat{\mu}(k)=\int_0^1 e^{-2\pi ikx}d\mu(x).
$$

\begin{theorem}\label{thm-main}
Assume that $E\subset[0,1]$ is a closed set which supports a probabilistic measure $\mu$
with the following properties:

\medskip

(A) $\mu([x,x+\epsilon])\leq C_1\epsilon^\alpha$ for all $0<\epsilon\leq 1$,

\medskip


(B) $|\widehat{\mu}(k)| \leq C_2 (1 - \alpha)^{-B} |k|^{-\frac{\beta}{2}}$ for all $k \ne 0$,  

\medskip\noindent
where $0<\alpha<1$ and $2/3<\beta\leq 1$.  If $\alpha>1-\epsilon_0$, where
$\epsilon_0>0$ is a sufficiently small constant depending only
on $C_1,C_2,B,\beta$, then $E$ contains a non-trivial 3-term arithmetic progression.

\end{theorem}

We note that if (A) holds with $\alpha=1$, then $\mu$
is absolutely continuous with respect to the Lebesgue measure, hence $E$ has positive
Lebesgue measure.  This case is already covered by the Lebesgue density argument
(see also Proposition \ref{preliminary-lemma} for a quantitative version). 

In practice, (B) will often be satisfied with $\beta$ very close
to $\alpha$.  It will be clear from the proof that the dependence
on $\beta$ can be dropped from the statement of the theorem if
$\beta$ is bounded from below away from $2/3$, e.g. $\beta>4/5$;
in such cases, the $\epsilon_0$ in Theorem \ref{thm-main} depends only on 
$C_1,C_2,B$.

The assumptions of Theorem \ref{thm-main} are in part suggested
by number-theoretic considerations, which we now describe briefly.
A theorem of Roth \cite{roth} 
states that if $A\subset\nn$ has {\em positive upper density}, i.e.
\begin{equation}\label{e-dens}
\overline{\lim}_{N\to\infty}\frac{\#(A\cap\{1,\dots,N\})}{N}>0,
\end{equation}
then $A$ must contain a non-trivial 
3-term arithmetic progression.  
Szemer\'edi's theorem extends this to $k$-term progressions.
It is well known that Roth's theorem fails without the assumption
(\ref{e-dens}), see \cite{salem-spencer}, \cite{behrend}.  However,
there are certain natural
cases when (\ref{e-dens}) may fail but the conclusion of Roth's theorem 
still holds. For example, there are variants of Roth's theorem for random
sets \cite{KLR}, \cite{tao-vu} and sets such as primes which resemble
random sets closely enough \cite{green-roth}, \cite{gt-2}.  
The key concept turns out to be {\em linear uniformity}. It is not hard
to prove that if the Fourier coefficients $\widehat{A}(k)$ of the characteristic function 
of $A$ are sufficiently small, depending on the size of $A$, then $A$ must contain 
3-term arithmetic progressions even if its asymptotic density is 0.  The 
Roth-type results mentioned above say that the same conclusion
holds under the weaker assumption that $A$ has an appropriate
{\em majorant} whose Fourier coefficients are sufficiently 
small (this is true for example if $A$ is a large subset of a random set).

If the universality of $A=\{0,1,2\}$ for sets of positive Lebesgue measure is viewed as
a continuous analogue of Roth's theorem, then its lower-dimensional analogue 
corresponds to Roth's theorem for integer sets of density 0 in $\nn$.
The above considerations suggest that such a result might hold under
appropriately chosen Fourier-analytic conditions on $E$ which could be interpreted
in terms of $E$ being ``random."  We propose Assumptions (A)-(B) of Theorem
\ref{thm-main} as such conditions.

To explain why Assumptions (A)-(B) are natural in this context, we give a brief
review of the pertinent background.
Let $\dim_H(E)$ denote the Hausdorff dimension of $E$.
Frostman's lemma (see e.g. \cite{falconer}, \cite{mattila-book}, \cite{wolff-lectures}) 
asserts that if $E\subset\rr$ is a compact set then
$$
\dim_H(E)=\sup\{\alpha\geq 0:\ \exists\hbox{ a probabilistic measure $\mu$ supported on $E$}$$
\vskip-8mm
$$
\hbox{ such that  (A) holds for some $C_1=C_1(\alpha)$}\}. 
$$
We also define the
{\em Fourier dimension} of $E\subset\rr$ as
$$
\dim_F(E)=\sup\{\beta\geq 0:\ \exists\hbox{ a probabilistic measure $\mu$ supported on $E$}$$
\vskip-8mm
$$
\hbox{ such that }|\widehat{\mu}(\xi)|\leq C
(1+|\xi|)^{-\beta/2}\hbox{ for all }\xi\in\rr \},
$$
where $\widehat{\mu}(\xi)=\int e^{-2\pi i\xi x}d\mu(x)$. 
Thus (A) implies that $E$ has Hausdorff dimension at least $\alpha$,
and (B) says that $E$ has Fourier dimension at least $2/3$.

It is known that 
\begin{equation}\label{i-e1}
\dim_F(E)\leq\dim_H(E)\hbox{ for all }E\subset \rr;
\end{equation}
in particular, a non-zero measure supported on $E$ cannot obey
(B) for any $\beta>\dim_H(E)$ (see (\ref{nu-e4})).  
It is quite common for the inequality in (\ref{i-e1}) to
be sharp: for instance, the middle-thirds Cantor set has Hausdorff
dimension $\log 2/\log 3$, but Fourier dimension 0.  Nonetheless, there are large
classes of sets such that
$$\dim_F(E)=\dim_H(E).$$
Such sets are usually called {\em Salem sets}.
It is quite difficult to construct explicit examples
of Salem sets with $0<\dim_H(E)<1$.  Such constructions are due
to Salem \cite{salem}, Kaufman \cite{kaufman}, Bluhm \cite{bluhm-1},
\cite{bluhm-2}; we give an alternative
construction in Section \ref{sec-examples}. On the other hand, 
Kahane \cite{kahane} showed that Salem sets are ubiquitous
among random sets, in the sense that images of compact sets under Brownian motion are
almost surely Salem sets.  This is one reason why we see Assumption (B) as a ``randomness" 
assumption and a good analogue of the number-theoretic concept of linear uniformity.

Assumptions (A)-(B) are closely related, but not quite equivalent,
to the statement that $E$ is a Salem set in the sense explained above.  
On the one hand, we do not have to
assume that the Hausdorff and Fourier dimensions of $E$ are actually
equal.  It suffices if (B) holds for {\em some} $\beta$, not necessarily
equal to $\alpha$ or arbitrarily close to it.  
On the other hand, we need to control the constants $C_1,C_2,B$,
as the range of $\alpha$ for which our theorem holds depends on these
constants.  (For example, we could set $B=0$, $C_1=C_2=100$, $\beta=4/5$; 
then our theorem states that if $\alpha$ is close enough to 1 (depending on 
the above choice of constants), then $\{0,1,2\}$ is universal for all sets which
support a measure $\mu$ obeying (A)-(B) with these values of $\alpha,\beta,B,C_1,C_2$.)

Thus we need to address the question of
whether measures obeying these modified assumptions can actually exist.
In Section \ref{sec-examples} we prove that given any $C_1>1$, $C_2>0$ and 
$0<\beta<\alpha<1$, there are subsets of $[0,1]$ which obey (A)-(B) with $B=0$ 
and with the given values of $C_1,C_2,\alpha,\beta$. Our construction is
based on probabilistic ideas similar to those of \cite{salem}, \cite{bluhm-1}, but 
simpler.  Salem's construction \cite{salem} does not produce explicit constants,
but we were able to modify his argument so as to show that, with large
probability, the examples in \cite{salem} obey (A)-(B) with $B=1/2$ and with 
$C_1,C_2$ independent of $\alpha$ for $\alpha$ close to 1.
Kahane's examples \cite{kahane} do not seem to obey (A)-(B) with
uniform constants; instead, they obey a condition similar to (A) but with an
additional logarithmic factor in $\epsilon^{-1}$, see e.g. \cite{DPRZ}.  We do not know whether 
our proof of Theorem \ref{thm-main} (specifically, the argument in Proposition
\ref{prop-restriction}) extends to this setting.  However, we are able
to give a more direct proof,
bypassing Proposition \ref{prop-restriction} and appealing
directly to Proposition \ref{prop-measure} instead, that Brownian image sets
do contain 3-term arithmetic progressions with positive probability bounded
from below uniformly in $\alpha$. 

Kaufman's set \cite{kaufman}, unlike those of Salem and Kahane, is completely
deterministic.  It is easily seen that if $x$ belongs to Kaufman's set, then so do
$2x,3x,\dots,$ in particular the set 
contains many $k$-term arithmetic progressions for any $k$.

The key feature of our proof is the use of a {\em restriction
estimate}.  Restriction estimates 
originated in Euclidean harmonic analysis, where they are
known for a variety of curved hypersurfaces 
(see e.g. \cite{stein-ha}).
In the paper \cite{mock} that inspired much of our work here, Mockenhaupt 
proved a restriction-type result for Salem sets in $\rr^d$.
Specifically, he proved that if $\mu$ obeys (A)-(B), then there is a 
restriction estimate of the form
\begin{equation}\label{e-restriction2}
\left[ \int |f|^2 d\mu \right]^{\frac{1}{2}} \leq A C_2^{\theta} C_1^{1 - \theta} 
||\widehat{f}||_{\ell^p(\mathbb Z)},
\end{equation}
for an appropriate range of $p$, analogous to the Stein-Tomas restriction 
theorem for the sphere in $\rr^n$ \cite{stein-beijing}, \cite{stein-ha},
\cite{tomas-1}, \cite{tomas-2}.  For our purposes, we will require
a variant of Mockenhaupt's estimate with uniform bounds on the constants,
which we derive in Section \ref{sec-restriction}.

While Mockenhaupt's work was primarily motivated by
considerations from Euclidean harmonic analysis, 
restriction estimates similar to (\ref{e-restriction2}) are also known
in number theory. Originally discovered by Bourgain \cite{bourg-89}, \cite{bourg-93},
they were recently used to prove Roth-type theorems for certain classes of integer 
sets of density zero.  Green \cite{green-roth} gave a proof of Roth's theorem in the
primes based on a restriction estimate for the primes.  Green's approach was refined
and extended further by Green-Tao \cite{gt-2}, and applied to a random set setting
by Tao-Vu \cite{tao-vu}.

Our proof of Theorem \ref{thm-main} extends the approach of \cite{green-roth},
\cite{gt-2}, \cite{tao-vu} to the continuous setting of sets of fractional dimension
for which a restriction estimate is available.  We rely on many of the
ideas from \cite{roth}, \cite{green-roth}, and particularly 
\cite{gt-2}, \cite{tao-vu}, such as the use of the trilinear form $\Lambda$
in a Fourier representation and a decomposition of the measure $\mu$ into
``random" and ``periodic" parts.  However, our actual argument is quite
different in its execution from those of \cite{gt-2} or \cite{tao-vu}. 
For instance, in \cite{green-roth}, \cite{gt-2}, \cite{tao-vu} the restriction
estimate is applied to the ``random" term, whereas
we use it to handle the ``periodic" part instead.

While Theorem \ref{thm-main} can be viewed as the analytic analogue of the
number-theoretic results just mentioned, it seems rather unlikely that 
our result could be deduced from them via a simple discretization procedure.
For instance, a $\delta$-neighbourhood of Keleti's set \cite{keleti}
contains many arithmetic progressions with common difference much 
greater than $\delta$; this immediately eliminates the simplest types of
discretization arguments.

Throughout the article, we
use $\#A$ to denote the cardinality of a finite set $A$, and $|E|$ to denote the 1-dimensional
Lebesgue measure of a set $E\subset\rr$. 

\bigskip

{\bf Acknowledgements.} 
We are grateful to Yuval Peres for bringing Erd\H{o}s's question to
our attention and for suggestions regarding Kahane's examples, and
to Mihalis Kolountzakis for further references including \cite{keleti}.
We also would like to thank Martin Barlow, Alex Iosevich, Nir Lev, 
Ed Perkins and Jim Wright for
helpful comments and suggestions.


\section{Outline of the proof of Theorem \ref{thm-main}}
\label{sec2}

Throughout this section we will assume that $\mu$ is a probabilistic measure supported
on a closed set $E\subset[1/3,2/3]$.   By scaling and translation, 
Theorem \ref{thm-main} extends to all closed $E\subset[0,1]$.
We will also invoke Assumptions (A) and (B) of
Theorem \ref{thm-main} where necessary.

We define the Fourier coefficients of a (possibly signed) measure $\sigma$ as
$$
\widehat{\sigma}(k)=\int_0^1 e^{-2\pi ikx}d\sigma(x).
$$
Given three signed measures $\mu_1$, $\mu_2$, $\mu_3$ on $[0,1]$, we define 
the trilinear form
$$
\Lambda(\mu_1,\mu_2,\mu_3)=\sum_{k=-\infty}^\infty \widehat{\mu_1}(k)
\widehat{\mu_2}(k)\widehat{\mu_3}(-2k).
$$
This notion is motivated by the Fourier-analytic proof of Roth's theorem
\cite{roth}, where a discrete version of $\Lambda$ is used to count the number of 
arithmetic progressions in a set of integers.  

We begin by considering the case of a measure 
$\mu$ absolutely continuous with respect to the Lebesgue measure.  Let $E\subset
[1/3,2/3]$ be a closed set, and let 
$\mu=df$, where $f$ is supported on $E$, $0\leq f\leq M$ and $\int_0^1 f(x)dx=1$. 
We will write
$$
\Lambda(\mu,\mu,\mu)=\Lambda(f,f,f)\hbox{ if }d\mu=f.
$$
\begin{lemma}\label{lambda-lemma1}
Let $f$ be a nonnegative bounded function supported on $[1/3,2/3]$. Then
\begin{equation}\label{mm-e1}
\Lambda(f,f,f)=\sum_{k=-\infty}^\infty \widehat{f}(k)^2\widehat{f}(-2k)
=2\int_0^1\int_0^1 f(x)f(y)f(\frac{x+y}{2})dxdy.
\end{equation}
\end{lemma}
\begin{proof}
This follows from a calculation identical to that in \cite{roth}.
Specifically, $\widehat{f}(k)^2$ is the Fourier transform of $f*f(x)
=\int_0^1 f(y)f(x-y)dy$ (note that since $E\subset[1/3,2/3]$, there is no need
to invoke addition mod 1 in the definition of $f*f$)
and $\widehat{f}(2k)$ is the Fourier transform
of $2f(x/2)$.  Hence by Plancherel's identity we have 
\begin{equation}\label{mm-e1a}
\begin{split}
\Lambda(f,f,f)&=\int 2f(\frac{x}{2})\int f(y)f(x-y)dydx\\
&=2\iint f(\frac{u+y}{2})\int f(y)f(u)dydu,
\end{split}
\end{equation}
where we changed variables $u=x-y$.
\end{proof}

The key result is the following proposition.
\begin{proposition}
Let $f : [0,1] \rightarrow [0,M]$ be a bounded function with 
\[ \int_{0}^{1} f(x) dx \geq \delta.  \] Then there exists $c = c(\delta,M)$ such that
\[ \Lambda_3(f,f,f) \geq c(\delta,M).  \] \label{preliminary-lemma}  
\end{proposition}

Proposition \ref{preliminary-lemma} is analogous to Varnavides's
theorem in number theory, a quantitative version of Roth's theorem
on 3-term arithmetic progressions.  It can be proved by following
exactly the proof of Varnavides's theorem as given e.g. in
\cite{tao-escorial}.

Proposition \ref{preliminary-lemma} implies 
in particular that the set 
$$
X=\{(x,y): f(x)>0,\ f(y)>0,\ f(\frac{x+y}{2})>0\}
$$
has positive 2-dimensional Lebesgue measure.  Since the 2-dimensional
measure of the line $x=y$ is 0, the set $X$ must contain many points $(x,y)$
with $x\neq y$.  Pick any such $(x,y)$, then the set $E$ contains
the non-trivial arithmetic progression
$\{x,\frac{x+y}{2},y\}$.  

We note that while the simple existence of 3-term arithmetic progressions
in sets of positive measure was an easy consequence of Lebesgue's density 
theorem, the quantitative result in Proposition \ref{preliminary-lemma}
is much more difficult and appears to require highly non-trivial methods
from number theory.

If the measure $\mu$ is singular, the infinite sum defining $\Lambda(\mu,\mu,\mu)$
does not have to converge in the first place.  Furthermore, there is no obvious analogue 
of (\ref{mm-e1}) and it is no longer clear how to interpret $\Lambda(\mu,\mu,\mu)$ in terms of
arithmetic progressions.   Nonetheless, if we assume that (B) holds for some $\beta>2/3$, then
\begin{equation}\label{nu-e3}
\sum_{k=-\infty}^\infty |\widehat{\mu}(-2k)\widehat{\mu}(k)^2|
\leq 1+C_2\sum_{k\neq 0}k^{-3\beta/2}<\infty.
\end{equation}
In particular, the sum defining $\Lambda(\mu,\mu,\mu)$ converges.
Moreover, we have the following.

\begin{proposition}\label{prop-measure}
Let $\mu$ be a probability measure supported on a compact set $E\subset[1/3,2/3]$ 
such that Assumption (B) holds for some $\beta\in(2/3,1]$.
Assume furthermore that $\Lambda(\mu,\mu,\mu)>0$.  Then there are $x,y\in E$ such that
$x\neq y$ and $\frac{x+y}{2}\in E$.

\end{proposition}

Theorem \ref{thm-main} now follows if we prove the next proposition.

\begin{proposition}\label{prop-main}
Let $\mu$ be a probability measure supported on a compact set $E\subset[0,1]$ 
such that Assumptions (A) and (B) hold for some $\alpha,\beta$ with
$0<\alpha<1$, $2/3<\beta<1$.
Then there are constants $0 < c_0, \epsilon_0 \ll 1$ (depending
only on $C_1,C_2, B, \beta$) 
such that if $1 - \epsilon_0\leq\alpha < 1$,  then
\[ \Lambda(\mu, \mu, \mu) \geq c_0. \] 
\end{proposition}

We prove Proposition \ref{prop-measure} in Section \ref{sec-measure}.
In Section \ref{sec-restriction} we prove the key restriction estimate 
needed in the proof of Proposition \ref{prop-main}; the latter 
follows in Section \ref{sec-mainproof}.  In the last three sections
we discuss three classes of examples of Salem sets: Salem's 
original construction, Kahane's Brownian images, and a new construction
due to the authors.


\section{Proof of Proposition \ref{prop-measure}}\label{sec-measure}

Let $\mu$ be as in Proposition \ref{prop-measure}.
We will prove the proposition by constructing a Borel measure $\nu$ on
$[0,1]^2$ such that
\begin{equation}\label{n-e1}
\nu([0,1]^2)>0,
\end{equation}
\begin{equation}\label{n-e2}
\nu\hbox{ is supported on the set }X=\{(x,y):\ x,y,\frac{x+y}{2}\in E\},
\end{equation}
\begin{equation}\label{n-e3}
\nu(\{(x,x):\ x\in[0,1]\})=0.
\end{equation}

In this section, it will be convenient to
work with the continuous Fourier transform 
$$
\widehat{\mu}(\xi)=\int_{-\infty}^\infty e^{-2\pi i \xi x}d\mu(x)
$$
instead of the Fourier series.  It is well known (see e.g. \cite{wolff-lectures}, Lemma 9.A.4)
that under the assumptions of Proposition \ref{prop-measure} we have
\begin{equation}\label{nu-e4}
|\widehat{\mu}(\xi)|\leq C'_2(1+|\xi|)^{-\beta/2},
\end{equation}
where $C'_2$ depends only on $C_2$.  In particular, since $\beta>2/3$, we have
$\widehat{\mu}\in L^3(\rr)$.
We also note that
\begin{equation}\label{nu-e4a}
\Lambda(\mu,\mu,\mu)=\int \whm^2(\xi)\whm(-2\xi)d\xi.
\end{equation}
Indeed, if $\mu$ is an absolutely continuous measure with density $f$,
then both sides of (\ref{nu-e4a}) are equal to 
$2\int_0^1\int_0^1 f(x)f(y)f(\frac{x+y}{2})dxdy$, by (\ref{mm-e1}) and 
the continuous analogue of the calculation in (\ref{mm-e1a}).  The general
case follows by a standard limiting argument.

Fix a non-negative Schwartz function $\psi$ on $\rr$ with 
$\int\psi=1$, let $\psieps=\epsilon^{-1}\psi(\epsilon^{-1}x)$,
and let $\phieps=\mu*\psieps$.  
Note that
\begin{equation}\label{nu-e5}
\widehat{\phieps}(\xi)=\widehat{\mu}(\xi)\widehat{\psieps}(\xi)
=\whm(\xi)\widehat{\psi}(\epsilon\xi).
\end{equation}
Since $\int\psi=1$, we have $\widehat{\psi}(0)=1$ and $|\widehat{\psi}(\xi)|\leq 1$
for all $\xi$.  Moreover, $\widehat{\psi}$ is continuous since $\psi$ is Schwartz.
It follows that
$\widehat{\psi}(\epsilon\xi)\to 1$ for all $\xi$ as $\epsilon\to 0$.  Hence
\begin{equation}\label{nu-e6}
|\widehat{\phieps}(\xi)|\leq\min(|\whm(\xi)|,|\widehat{\psieps}(\xi)|),
\end{equation}
\begin{equation}\label{nu-e7}
\widehat{\phieps}(\xi)\to \whm(\xi)\hbox{ pointwise as }\epsilon\to 0.
\end{equation}

We define a linear functional $\nu$ acting on
functions $f:[0,1]^2\to\rr$ by the formula
\begin{equation}\label{nu-e1}
\langle \nu,f\rangle =\lim_{\epsilon\to 0}\iint f(x,y)\phieps(\frac{x+y}{2})d\mu(x)d\mu(y).
\end{equation}
Clearly, $\langle \nu,f\rangle\geq 0$ if $f\geq 0$.

\begin{lemma}\label{nu-lemma1}
The limit in (\ref{nu-e1}) exists for all continuous functions $f$ on $[0,1]^2$.
Moreover,
\begin{equation}\label{nu-e8}
|\langle \nu,f\rangle|\leq C\|f\|_\infty,
\end{equation}
where $C$ depends on $\mu$ but is independent of $f$.
\end{lemma}

\begin{proof}
Suppose that $f$ is continuous and bounded by $M$ on $[0,1]^2$.  Then, by
the same calculation as in (\ref{mm-e1a}),
\begin{equation}\label{nu-e2}
\begin{split}
\iint \Big|f(x,y)\phieps(\frac{x+y}{2})\Big|d\mu(x)d\mu(y)
 &\leq M\iint \phieps(\frac{x+y}{2})d\mu(x)d\mu(y)\\
 &=\frac{M}{2}\int \widehat{\phieps}(-2\xi)\widehat{\mu}(\xi)^2d\xi\\
 &\leq \frac{M}{2}\int |\widehat{\mu}(-2\xi)\widehat{\mu}(\xi)|^2d\xi<CM/2,
\end{split}
\end{equation}
where at the last step we used (\ref{nu-e6}) and (\ref{nu-e4}).  This implies that
if the limit in (\ref{nu-e1}) exists, then (\ref{nu-e8}) holds.  

It remains to prove the existence of the limit.
The Schwartz functions $f$ with $\widehat{f}\in C_0^\infty(\rr^2)$ are dense in
$C([0,1]^2)$ in the $L^\infty$ norm.  If we prove that the limit in (\ref{nu-e1}) exists
for such functions, it will follow from (\ref{nu-e8}) that it also exists for all 
continuous functions on $[0,1]^2$.

Let  $f$ be a Schwartz function on $\rr^2$ with $\supp\widehat{f}\subset\{(\eta_1,\eta_2):
0\leq |\eta_1|+|\eta_2|\leq R\}$.
By Plancherel's identity, we have
\begin{equation*}
\begin{split}
\langle \nu,f\rangle &=\lim_{\epsilon\to 0}\iint f(x,y)\phieps(\frac{x+y}{2})d\mu(x)d\mu(y)\\
&=\lim_{\epsilon\to 0}\iint \Phi_\epsilon(\eta_1,\eta_2)\whm(\eta_1)\whm(\eta_2)d\eta_1d\eta_2,
\end{split}
\end{equation*}
where
\begin{equation*}
\begin{split}
\Phi_\epsilon(\eta_1,\eta_2)
&=\iint f(x,y)\phieps(\frac{x+y}{2})e^{-2\pi i(x\eta_1+y\eta_2)}dxdy\\
&=\iiint f(x,y)\widehat{\phieps}(\xi)e^{-2\pi i(x\eta_1+y\eta_2-\frac{x+y}{2}\xi)}d\xi dxdy\\
&=\int \widehat{\phieps}(\xi)\Big(\iint f(x,y)e^{-2\pi i(x(\eta_1-\xi/2)+y(\eta_2-\xi/2)}dxdy\Big)d\xi\\
&=\int \widehat{\phieps}(\xi)\widehat{f}\Big(\eta_1-\frac{\xi}{2},\eta_2-\frac{\xi}{2}\Big)d\xi.
\end{split}
\end{equation*}
Hence
\begin{equation}\label{nu-e10}
\langle \nu,f\rangle =\lim_{\epsilon\to 0}\iiint 
 \widehat{\phieps}(\xi)\widehat{f}\Big(\eta_1-\frac{\xi}{2},\eta_2-\frac{\xi}{2}\Big)
 \whm(\eta_1)\whm(\eta_2)d\eta_1d\eta_2d\xi.
\end{equation}
By (\ref{nu-e6}), the integrand in (\ref{nu-e10}) is bounded in absolute value by
\begin{equation}\label{nu-e11}
\Big|\widehat{\mu}(\xi)\widehat{f}\Big(\eta_1-\frac{\xi}{2},\eta_2-\frac{\xi}{2}\Big)
 \whm(\eta_1)\whm(\eta_2)\Big|
\end{equation}
for all $\epsilon>0$. 

We claim that
\begin{equation}\label{nu-e12}
\iiint \Big|\widehat{\mu}(\xi)\widehat{f}\Big(\eta_1-\frac{\xi}{2},\eta_2-\frac{\xi}{2}\Big)
 \whm(\eta_1)\whm(\eta_2)\Big|d\xi d\eta_1 d\eta_2<\infty.
\end{equation}
Assuming (\ref{nu-e12}), it follows from (\ref{nu-e7}) and the dominated convergence theorem
that the limit in (\ref{nu-e10}) exists and is equal to
\begin{equation*}
\iiint \widehat{\mu}(\xi)\widehat{f}\Big(\eta_1-\frac{\xi}{2},\eta_2-\frac{\xi}{2}\Big)
 \whm(\eta_1)\whm(\eta_2)d\xi d\eta_1 d\eta_2,
\end{equation*}
which proves the lemma.

We now prove the claim (\ref{nu-e12}).  Note first that by the support properties of $\widehat{f}$,
the integral in (\ref{nu-e12}) is in fact taken over the set
$$\Omega=\Big\{(\xi,\eta_1,\eta_2):\ \Big|\eta_1-\frac{\xi}{2}\Big|\leq R,\ 
\Big|\eta_2-\frac{\xi}{2}\Big|\leq R\Big\}.
$$
Let $1< p,p'<\infty$ be exponents such that $\frac{1}{p}+\frac{1}{p'}=1$ and
\begin{equation}\label{nu-e14}
p'\beta/2>1.
\end{equation}
Then $\whm(\xi)\in L^{p'}(\rr,d\xi)$ and, by H\"older's inequality, the left side of 
(\ref{nu-e12}) is bounded by
\begin{equation}\label{nu-e20}
\|\whm\|_{L^{p'}(d\xi)}\cdot\|F(\xi)\|_{L^p(d\xi)},
\end{equation}
where
\begin{equation}\label{nu-e15}
F(\xi)=\iint_{|\eta_1-\frac{\xi}{2}|\leq R,\ 
|\eta_2-\frac{\xi}{2}|\leq R}
\Big|\widehat{f}\Big(\eta_1-\frac{\xi}{2},\eta_2-\frac{\xi}{2}\Big)
 \whm(\eta_1)\whm(\eta_2)\Big|d\eta_1 d\eta_2.
\end{equation}
By H\"older's inequality, we have
$$
F(\xi)\leq \|\widehat{f}\|_{L^{q'}(d\xi)}
\Big(\iint |\whm(\eta_1)\whm(\eta_2)|^q d\eta_1 d\eta_2\Big)^{1/q},
$$
where $1< q,q'<\infty$ are exponents such that $\frac{1}{q}+\frac{1}{q'}=1$, and
the double integral is taken over the same region as in (\ref{nu-e15}).
On that region we have $|\eta_1-\eta_2|\leq 2R$ and
$$
|\whm(\eta_1)|\leq C(1+|\eta_1|)^{-\beta/2}\leq C_R (1+|\xi|)^{-\beta/2}
\leq C_R^2 (1+|\eta_2|)^{-\beta/2},
$$
and similarly with indices 1 and 2 interchanged.  Hence
\begin{equation*}
\begin{split}
F(\xi)&\leq C(f,q)\Big(\int_{-\infty}^\infty (1+|\eta_2|)^{-2}
\int_{\eta_2-2R}^{\eta_2+2R} (1+|\xi|)^{2-\beta q}d\eta_1d\eta_2\Big)^{1/q}\\
& \leq C'(f,q) (1+|\xi|)^{\frac{2}{q}-\beta}.\\
\end{split}
\end{equation*}
By (\ref{nu-e20}), it suffices to find exponents $p,q$ such that (\ref{nu-e14}) holds
and $F(\xi)\in L^p$.  If
\begin{equation}\label{nu-e21}
\beta p>1
\end{equation}
and if $q$ is chosen large enough, then $(\beta-\frac{2}{q})p>1$ and $F(\xi)$ is
$L^p$-integrable as required.  Finally, if
$$
1-\frac{\beta}{2}<\frac{1}{p}<\beta,
$$
which is possible whenever $\beta>2/3$, both (\ref{nu-e21}) and (\ref{nu-e14}) hold.
This completes the proof of the lemma.

\end{proof}

By Lemma \ref{nu-lemma1} and the Riesz representation theorem, (\ref{nu-e1})
defines a measure on $[0,1]^2$.  We will now prove that $\nu$ has the desired properties
(\ref{n-e1})--(\ref{n-e3}).

\bigskip

\noindent{\it Proof of (\ref{n-e1}).}  We write
 \begin{equation*}
\begin{split}
\langle \nu,1\rangle &=
\lim_{\epsilon\to 0} \iint \phieps(\frac{x+y}{2})d\mu(x)d\mu(y)\\
 &=\lim_{\epsilon\to 0}\int \widehat{\phieps}(-2\xi)\widehat{\mu}(\xi)^2d\xi\\
 &=\int \widehat{\mu}(-2\xi)\widehat{\mu}(\xi)^2d\xi,
\end{split}
\end{equation*}
by (\ref{nu-e6}), (\ref{nu-e7}), (\ref{nu-e4}) and the dominated convergence
theorem.  By (\ref{nu-e4a}), the last integral is equal to
$\Lambda(\mu,\mu,\mu)$, hence is positive
as claimed.

\bigskip
\noindent{\it Proof of (\ref{n-e2}).} Since $E$ is closed, $X$ is closed.  Let $f$ be 
a continuous function with $\supp f$ disjoint from $X$, then $\dist(\supp f,X)>0$. 
Using a partition of unity, we write $f=\sum f_j$, where $f_j$ are continuous and
for each $j$ at least one of the following holds:
$$\dist(\supp f_j,\ E\times[0,1])>0,$$
$$\dist(\supp f_j,\ [0,1]\times E)>0,$$
$$\dist\Big(\Big\{\frac{x+y}{2}:\ (x,y)\in\supp f_j\Big\},\ E\Big)>0.$$
It suffices to prove that $\langle \nu,f_j\rangle=0$ for all $j$.
In the first two cases, we have $\mu\times\mu(\supp f_j)=0$, hence the integral
in (\ref{nu-e1}) is 0 for each $\epsilon$.  In the last case, we have $\phieps(\frac{x+y}{2})
\to 0$ pointwise on the support of $f_j$, hence by the dominated convergence theorem 
the integral in (\ref{nu-e1}) converges to 0 as $\epsilon\to 0$, as required.

\bigskip
\noindent{\it Proof of (\ref{n-e1}).}  
It suffices to prove that 
\begin{equation}\label{nu-e100}
\nu(\{(x,y): \ |x-y|\leq\delta\})\to 0\hbox{ as }\delta\to 0.
\end{equation}
Let $\chi:\rr\to[0,\infty)$ be a Schwartz function such that $\chi\geq 0$, 
$\chi\geq 1$ on $[-1,1]$
and $\supp\widehat{\chi}\subset\{|\xi|\leq R\}.$  We will prove that
\begin{equation}\label{nu-e101}
\langle \nu,\chi_\delta(x-y)\rangle \to 0\hbox{ as }\delta\to 0,
\end{equation}
where
$\chi_\delta(t)=\chi(\delta^{-1}t)$. This will clearly imply (\ref{nu-e100}).

We write
\begin{equation*}
\begin{split}
&\langle \nu,\chi_\delta(x-y)\rangle \\
&=\lim_{\epsilon\to 0}\iint \chi_\delta(x-y)\phieps(\frac{x+y}{2})d\mu(x)d\mu(y)\\
&=\lim_{\epsilon\to 0}\iint\Big(\int \widehat{\chi_\delta}(\xi)e^{2\pi i(x-y)\xi}d\xi\Big)
\Big(\int\widehat{\phieps}(\eta)e^{\pi i(x+y)\eta}d\eta\Big)d\mu(x)d\mu(y)\\
&=\lim_{\epsilon\to 0}\iint \widehat{\chi_\delta}(\xi)\widehat{\phieps}(\eta)
\whm(\xi-\frac{\eta}{2})\whm(-\xi-\frac{\eta}{2})d\xi d\eta.
\end{split}
\end{equation*}
We claim that 
\begin{equation}\label{nu-e103}
 \iint \Big|\widehat{\chi_\delta}(\xi)\whm(\eta)
\whm(\xi-\frac{\eta}{2})\whm(-\xi-\frac{\eta}{2})\Big|d\xi d\eta
\leq C\delta^s
\end{equation}
for some $s>0$.  Here and throughout the proof, we will use $C,C',\dots$ to denote 
positive constants which may change from line to line and may depend on $\mu$,
$\chi$ and $R$, but are always uniform in $\delta$. 
Assuming (\ref{nu-e103}), we have
$$
\langle \nu,\chi_\delta(x-y)\rangle
 =\iint \widehat{\chi_\delta}(\xi)\whm(\eta)
\whm(\xi-\frac{\eta}{2})\whm(-\xi-\frac{\eta}{2})d\xi d\eta
$$
by (\ref{nu-e7}) and the dominated convergence theorem.  By (\ref{nu-e103}) again,
the last integral is bounded by $C\delta^s$, hence (\ref{nu-e101}) follows.

It remains to prove (\ref{nu-e103}).  Since $\widehat{\xi_\delta}(\xi)=\delta
\widehat{\chi}(\delta\xi)$ is supported on $|\xi|\leq \delta^{-1}R$, we can bound
the integral in (\ref{nu-e103}) by
\begin{equation}\label{nu-e104}
\begin{split}
&C\delta  \iint_{|\xi|\leq\delta^{-1}R} \Big|\whm(\eta)
\whm(\xi-\frac{\eta}{2})\whm(-\xi-\frac{\eta}{2})\Big|d\xi d\eta\\
&=C\delta  \int |\whm(\eta)|\Big(\int_{|\xi|\leq\delta^{-1}R} \Big|
\whm(\xi-\frac{\eta}{2})\whm(-\xi-\frac{\eta}{2})\Big|d\xi\Big) d\eta\\
&\leq C\delta  \|\whm\|_{p'}\|F(\eta)\|_{L^{p}(d\eta)},
\end{split}
\end{equation}
where $1<p,p'<\infty$ are exponents such that $\frac{1}{p}+\frac{1}{p'}=1$, and
$$
F(\eta)=\int_{|\xi|\leq\delta^{-1}R} \Big|
\whm(\xi-\frac{\eta}{2})\whm(-\xi-\frac{\eta}{2})\Big|d\xi.
$$
By (\ref{nu-e4}), $\|\whm\|_{p'}$ is finite whenever
\begin{equation}\label{nu-e106}
p'\beta/2>1.
\end{equation}
We now estimate $F(\eta)$.  Let $1<q,q'<\infty$ are exponents such that 
$\frac{1}{q}+\frac{1}{q'}=1$, then by H\"older's inequality and (\ref{nu-e4}) again we have
\begin{equation}\label{nu-e107}
\begin{split}
F(\eta)&\leq
\Big(\int_{|\xi|\leq\delta^{-1}R} 1d\xi\Big)^{1/q'}
\Big(\int_{|\xi|\leq\delta^{-1}R} \Big|
\whm(\xi-\frac{\eta}{2})\whm(-\xi-\frac{\eta}{2})\Big|^q d\xi\Big)^{1/q}\\
&\leq C\delta^{-1/q'},
\end{split}
\end{equation}
assuming that
\begin{equation}\label{nu-e108}
q\beta>1.
\end{equation}
Moreover, on the region $|\eta|>C'\delta^{-1}$ for sufficiently large $C'$ we have
$$\Big|\xi\pm\frac{\eta}{2}\Big|\sim\Big|\frac{\eta}{2}\Big|\hbox{ for all }
|\xi|\leq \delta^{-1}R,$$
hence the calculation in (\ref{nu-e107}) yields the stronger estimate
\begin{equation}\label{nu-e109}
\begin{split}
F(\eta)&\leq C\delta^{-1/q'}
\Big(\int_{|\xi|\leq\delta^{-1}R} (1+\eta)^{-\beta q} d\xi\Big)^{1/q'}\\
&\leq C\delta^{-1/q'}
\Big(\delta^{-1}(1+\eta)^{-\beta q} \Big)^{1/q'}\\
&\leq C\delta^{-1}|\eta|^{-\beta}.
\end{split}
\end{equation}
Now we can bound the $L^p$-norm of $F$:
$$
\Big(\int F(\eta)^pd\eta\Big)^{1/p}
\leq\Big(\int_{|\eta|\leq C'\delta^{-1}}F(\eta)^p d\eta\Big)^{1/p}
+\Big(\int_{|\eta|\geq C'\delta^{-1}}F(\eta)^p d\eta\Big)^{1/p},
$$
$$
\Big(\int_{|\eta|\leq C'\delta^{-1}}F(\eta)^p d\eta\Big)^{1/p}
\leq C\Big(\delta^{-1}(\delta^{-1/q'})^p\Big)^{1/p}=C\delta^{-\frac{1}{p}-\frac{1}{q'}},
$$
\begin{equation*}
\begin{split}
\Big(\int_{|\eta|\geq C'\delta^{-1}}F(\eta)^p d\eta\Big)^{1/p}
&= C\Big(\int_{C'\delta^{-1}}^\infty
\delta^{-p}|\eta|^{-\beta p} d\eta\Big)^{1/p}\\
&= C\delta^{-1}\Big((\delta^{-1})^{-\beta p+1}\Big)^{1/p}\\
&=C\delta^{-1+\beta-\frac{1}{p}}.
\end{split}
\end{equation*}
Returning now to (\ref{nu-e104}), we see that the integral in (\ref{nu-e103}) 
is bounded by
\begin{equation}\label{nu-e110}
C(\delta^{1-\frac{1}{p}-\frac{1}{q'}}+\delta^{\beta-\frac{1}{p}}).
\end{equation}
The exponent of $\delta$ in (\ref{nu-e110}) is positive if 
\begin{equation}\label{nu-e111}
\frac{1}{p}+\frac{1}{q'}\leq 1\hbox{ and }\beta >\frac{1}{p}.
\end{equation}
A short calculation shows that we can find $p,q$ obeying all of (\ref{nu-e111}),
(\ref{nu-e106}), (\ref{nu-e107}) whenever $\beta>2/3$.  This ends the proof
of (\ref{nu-e104}) and (\ref{n-e3}).


\section{A restriction estimate}\label{sec-restriction}
The key part of the proof of Proposition \ref{prop-main} is a Fourier restriction estimate 
for measures $\mu$ satisfying Assumptions (A) and (B)
of Theorem \ref{thm-main}.
Such estimates were first obtained by Mockenhaupt in \cite{mock}, following Tomas's
proof of a restriction estimate for the sphere in $\mathbb R^n$ \cite{tomas-1}. The main 
result in \cite{mock} is an $L^p(dx) \rightarrow L^2(d \mu)$ bound for the Fourier 
restriction operator in the range $1 \leq p < 2(2- 2 \alpha + \beta)/(4(1 - \alpha) + \beta)$. 
For our application, we are interested in values of $p$ near the endpoint. However, it 
is necessary for us to keep track of the behavior of the operator norm (in terms of 
$\alpha$ and $\beta$) near the endpoint, in particular ensuring that it stays bounded 
for $\alpha, \beta$ close to 1, and explicitly deriving its dependence on $C_1$ and 
$C_2$. Mockenhaupt's analysis yields a bound for the operator norm that blows up 
at the endpoint. Our result in contrast gives a uniform bound for the operator norm, 
though not ``at'' the endpoint, only ``near'' it.

\begin{proposition}\label{prop-restriction}
Let $\frac{2}{3}< \alpha, \beta \leq 1$, and let $\mu$ be a probability measure 
supported on $E \subset [0, 1]$ and obeying Assumptions (A) and (B) of Theorem
\ref{thm-main}.
Then there exists an absolute constant $A$ such that
\begin{equation}\label{e-restriction}
\left[ \int |f|^2 d\mu \right]^{\frac{1}{2}} \leq 2^{6B}AC_2^{\theta} C_1^{1 - \theta} 
||\widehat{f}||_{\ell^p(\mathbb Z)},
\end{equation}
where $p = \frac{2 (\beta + 4(1 - \alpha))}{\beta + 8(1 - \alpha)}$  and $\theta = \frac{2}{p}-1$. 
\end{proposition}

\begin{proof}
By duality, it suffices to show that the operator $T$ defined by
$Tg = g \ast \widehat{\mu}$ maps $\ell^p(\mathbb Z) \rightarrow \ell^{p'}(\mathbb Z)$,
with operator norm bounded by $AC_2^{\theta}C_1^{1 - \theta}$. We prove this
by complex interpolation. Following an approach analogous to Stein's proof
of the endpoint restriction estimate in the
Stein-Tomas theorem \cite{stein-beijing}, \cite{stein-ha}, 
we embed $T$ in the family of operators $T_sg = g \ast \widehat{K}_s$,
where $\widehat{K}_s(n) = \zeta_s(n) \widehat{\mu}(n)$. Here 
\[ \zeta_s(\xi) = e^{s^2}\frac{3(1 - \alpha)s}{3(1 - \alpha) - s}
\int e^{-2 \pi i x \xi} |x|^{-1 + s} \eta(x) \, dx, \]  
with $\eta \in C_0^{\infty}(\mathbb R)$ supported in $[-1,1]$, $0 \leq \eta \leq 1$,
and  $\eta \equiv 1$ near the origin. Then $K_s$, initially defined on the half-plane
Re$(s) > 0$ admits an analytic continuation in the interior of the strip
$-\frac{\beta}{2} \leq \text{Re}(s) \leq 2(1- \alpha)$ and is continuous up to
its boundary. Moreover, the arguments in \cite{stein-ha}(page 381-382) show that
$\zeta_0(x) \equiv 1$. Thus the desired result will follow from the
three-lines lemma if we establish that $||T_s||_{\ell^1 \rightarrow \ell^{\infty}}
\leq AC_2$ for Re$(s) = -\frac{\beta}{2}$, and
$||T_s||_{\ell^2 \rightarrow \ell^2} \leq AC_1$ for Re$(s) = 2(1- \alpha)$.
Specifically, we need the two estimates 
\begin{equation}\label{res-e1}
\sup_n |\widehat{K}_s(n)| \leq AC_2(1-\alpha)^{-B} \quad \text{ for } \quad \text{Re}(s) = -\frac{\beta}{2},  
\end{equation}
\begin{equation}\label{res-e2}
\sup_x |K_s(x)| \leq AC_1 \quad \text{ for } \quad \text{Re}(s) = 2(1-\alpha). 
\end{equation}   

For Re$(s) = -\frac{\beta}{2}$, we have from \cite{stein-ha}   
\[ |\zeta_s(\xi)| \leq A (1 + |\xi|)^{\frac{\beta}{2}}, \]
for some constant $A$ uniform in $\alpha$, $\beta$ for $\frac{2}{3} \leq \alpha, \beta \leq 1$. 
This together with (B) yields (\ref{res-e1}). 

For Re$(s) = 2(1- \alpha)$, we use Assumption (A) to obtain
{\allowdisplaybreaks \begin{align*}
|K_s(x)| &= |\widehat{\zeta}_s \ast \mu(x)|\\
 &= \left|e^{s^2} \frac{3(1-\alpha)s}{3(1 - \alpha) - s} \right| 
 \left| \int |x - x'|^{-1 + s} \eta(x-x') d\mu(x') \right| \\ 
&\leq \left|e^{s^2} \frac{3(1-\alpha)s}{3(1 - \alpha) - s} \right| 
\sum_{j \geq 0} \left| \int_{|x - x'| \sim 2^{-j}} |x - x'|^{-1 + s} \eta(x-x') d\mu(x') \right| \\ 
&\leq \left|e^{s^2} \frac{3(1-\alpha)s}{3(1 - \alpha) - s} \right| 
\sum_{j \geq 0} 2^{-j(-1 + 2(1 - \alpha))} C_1 2^{-j\alpha} \\ 
&\leq C_1 \left|e^{s^2} \frac{s}{3(1 - \alpha) - s} \right| (1 - \alpha) \sum_{j \geq 0} 2^{-j(1 - \alpha)} \\ 
&\leq C_1 \sup_{t \in \mathbb R} \left[ e^{4(1 - \alpha)^2 - t^2} 
\frac{2(1 - \alpha) + |t|}{(1 - \alpha)+|t|} \right] \frac{(1-\alpha)}{1 - 2^{-(1 - \alpha)}} \\ &\leq AC_1.
\end{align*}}
Hence 
\begin{align*}
\|T_s\|_{\ell^p\rightarrow\ell^{p'}}
&\leq (AC_2(1-\alpha)^{-B})^{\theta}\,(AC_1)^{1-\theta}\\
&=AC_1^{1-\theta} C_2^{\theta}(1-\alpha)^{\frac{-4B(1-\alpha)}{4(1-\alpha)+\beta}}\\
&\leq AC_1^{1-\theta} C_2^{\theta}2^{6B}.\\
\end{align*}
At the last step we used that $(1-\alpha)^{-(1-\alpha)}\leq 2$ for all $2/3<\alpha<1$, hence the last factor is bounded by $2^{4B/(4-4\alpha+\beta)}\leq 2^{4B/\beta}\leq 2^{6B}$.
This proves (\ref{res-e2}) and completes the proof of the proposition.
\end{proof}

\section{Proof of Proposition \ref{prop-main}}
\label{sec-mainproof}

The proof of the proposition follows roughly the scheme in \cite{green-roth},
\cite{gt-2}, \cite{tao-vu}. We will find a decomposition 
$\mu=\mu_1+\mu_2$, where $\mu_1$ is absolutely continuous with
bounded density, and $\mu_2$ is irregular but ``random" in the sense that
it has very small Fourier coefficients. 
We then write
$$\Lambda(\mu,\mu,\mu)=\sum_{i,j,k=1}^2 \Lambda(\mu_i,\mu_j,\mu_k).$$
The main contribution will come from the term 
$\Lambda(\mu_1,\mu_1,\mu_1)$, which we will bound from below using
Proposition \ref{preliminary-lemma}.  The remaining terms will be
treated as error terms and will be shown to be small compared to
$\Lambda(\mu_1,\mu_1,\mu_1)$.

We begin by defining $\mu_1$ and $\mu_2$.  Let $K_N$ denote the Fej\'er kernel on 
$[0,1]$, namely 
\begin{equation} 
K_{N}(x) = \sum_{n = -N}^N \left(1 - \frac{|n|}{N+1} \right) e^{2\pi inx} = 
\frac{1}{N+1} \frac{\sin^2 \left((N+1)\pi x \right)}{\sin^2(\pi x)}. 
\label{fejer1} \end{equation} 
It follows easily from (\ref{fejer1}) that $K_N \geq 0$ and $\int_0^1 K_N = 1$ for every 
$N \geq 1$, $N \in \mathbb N$. Moreover, 
\begin{equation} K_{2N}(x) = \frac{1}{2N+1} D_N^2(x), \label{fejer2} \end{equation}
where 
$$D_N (x) = \sum_{|n| \leq N} e^{2\pi inx}
= \frac{\sin((2N +1)\pi x)}{\sin(\pi x)}$$ 
is the Dirichlet kernel.

Let $N \gg 1$ be a large constant to be determined shortly.  Let 
$$\mu_1(x) = K_{2N} \ast \mu(x),$$
and write $\mu = \mu_1 + \mu_2$. Clearly, $\mu_1 \geq 0$ and $\int_0^1 \mu_1(x) dx = 1$. 
We claim that if $\alpha,\beta$ are close enough to 1, then $N$ can be chosen 
so that
\begin{equation}\label{xx-1}
0\leq \mu_1(x) \leq M
\end{equation}
for some fixed constant $M$ depending on $C_1$ but 
{\em{independent of }} $N$. Indeed, by (\ref{fejer2}) and Proposition \ref{prop-restriction}
there is an absolute constant $A$ such that 
\begin{align*} \mu_1(x) = \int K_N(x-y) d\mu(y) 
&= \frac{1}{2N+1} \int |D_N(x-y)|^2 d\mu(y) \\  
&\leq \frac{2^{6B}A C_2^{\theta} C_1^{1 - \theta}}{2N+1}||\widehat{D}_N||_{\ell^p}^2 \\ 
& \leq 2^{6B}AC_1^{2/p'} (C_2(2N+1))^{\frac{2}{p} - 1},  \end{align*} 
where
$\displaystyle{p = \frac{2 (\beta + 4(1 - \alpha))}{\beta + 8(1 - \alpha)}}$. 

We may assume that $C_1,C_2\geq 1$, and we continue to assume that $\beta>2/3$.
For $\alpha \geq 1- (\ln C_2)^{-1}$, set $N =\lfloor C_2^{-1}e^{\frac{1}{1 - \alpha}}\rfloor$.  
Then $N\geq 1$ and $N\to\infty$ as $\alpha\to 1$.
Since $\displaystyle{\frac{2}{p}-1 = \frac{4(1 - \alpha)}{\beta + 4(1 - \alpha)}}$
and $\displaystyle{\frac{2}{p'} = \frac{2\beta}{\beta + 4(1 - \alpha)}\leq 1}$,
\begin{equation*}
\begin{split}
\mu_1(x)&\leq 2^{6B+2}A C_1^{2/p'} (C_2N)^{\frac{2}{p}-1}
\leq 2^{6B+2}A C_1
e^{\frac{1}{1-\alpha}\cdot \frac{4(1-\alpha)}{\beta+4(1-\alpha)}}\\
&\leq 2^{6B+2}AC_1e^6.
\end{split}
\end{equation*}
By Proposition \ref{preliminary-lemma}, 
there is a constant $c_0> 0$ (depending on $A,B,C_1$) 
such that $\Lambda_3(\mu_1, \mu_1, \mu_1) \geq 2c_0$.   

It remains to verify that the error terms are negligible: 
if $\beta$ is sufficiently close to 1, and at least one of the indices $i_1, i_2, i_3$ equals 2, then $$\Lambda(\mu_{i_1}, \mu_{i_2}, \mu_{i_3}) \leq c_0/8.$$
We only consider the cases $\Lambda(\mu_1, \mu_1, \mu_2)$ and 
$\Lambda(\mu_2, \mu_2, \mu_2)$; the other cases are similar and left to the reader. 

By the definition of $\mu_2$, we have
$\widehat{\mu}_2(0) = 0$ and
\[\widehat{\mu}_2(n) = \min \left(1, \frac{|n|}{2N+1} \right) \widehat{\mu}(n).
\] 
Hence
\begin{equation}\label{dec-e1}
\begin{split}
|\Lambda(\mu_1, \mu_1, \mu_2)| 
&\leq \sum_{0 \leq |n| \leq 2N} |\widehat{\mu}_1(n)|^2 |\widehat{\mu}_2(-2n)| \\ 
&\leq C_2^3(1-\alpha)^{-3B} \sum_{0<|n| \leq 2N} |n|^{-\beta} \frac{|n|}{2N+1}(2|n|)^{- \frac{\beta}{2}} \\
&=\frac{C_2^3(1-\alpha)^{-3B}}{2N+1} \sum_{0<|n| \leq 2N} |n|^{1-\frac{3\beta}{2}}
\\
&\leq 4C_2^3(1-\alpha)^{-3B}N^{1-\frac{3\beta}{2}}
\\
&\leq 4C_2^3(1-\alpha)^{-3B}
(C_2^{-1}e^{\frac{1}{1-\alpha}}-1)^{1-\frac{3\beta}{2}},\\
\end{split}
\end{equation}
which can be made arbitrarily small as $\alpha\to 1$ and $N\to\infty$.
Similarly, 
\begin{equation}\label{dec-e2}
\begin{split}
&|\Lambda(\mu_2, \mu_2, \mu_2)|\\ 
&\leq \frac{C_2^3(1-\alpha)^{-3B}}{( 2N+1)^{3}}
\sum_{|n| \leq 2N } |n|^{ 3 - \frac{3\beta}{2}} + 
C_2^3(1-\alpha)^{-3B} \sum_{|n| \geq 2N+1}|n|^{-\frac{3 \beta}{2}} \\ 
&\leq \frac{3\beta+2}{3\beta-2}C_2^3(1-\alpha)^{-3B}
N^{1-\frac{3\beta}{2}}\\
&\leq \frac{3\beta+2}{3\beta-2}C_2^3(1-\alpha)^{-3B}
(C_2^{-1}e^{\frac{1}{1-\alpha}}-1)^{1-\frac{3\beta}{2}}
\rightarrow 0 \text{ as } N \rightarrow \infty. 
\end{split}
\end{equation}


\section{A construction of Salem-type sets}
\label{sec-examples}

We now give a probabilistic construction of examples of sets 
equipped with natural measures obeying Assumptions
(A)-(B) of Theorem \ref{thm-main}.  The idea is reasonably simple. 
Start with the interval $[0,1]$ equipped with the Lebesgue measure. 
Subdivide it into $M_1=KN_1$ intervals of equal length, then choose
$Kt_1$ of them at random and assign weight $(Kt_1)^{-1}$ to each one. 
For a {\em generic} choice of subintervals, the Fourier transform of
the resulting density function is close to the Fourier transform
of ${\bf 1}_{[0,1]}$.  Now subdivide each of the intervals chosen
at the first step into $N_2$ subintervals of equal length, choose $t_2$ of
them at random, assign weight $(Kt_1t_2)^{-1}$ to each one.  Continue to 
iterate the procedure, taking care at each step to keep the Fourier
transform of the new density function as close as possible to the
previous one.
In the limit, we get a Cantor-type set equipped 
with a natural measure  $\mu$ such that $\widehat{\mu}$ has
the required decay. 

In the construction below, we will let $N_1=N_2=\dots
=N$ and $t_1=t_2=\dots=t$ be fixed.  This will produce Salem sets of
dimension $\alpha=\frac{\log t}{\log N}$.  
It should be possible to use the same
argument to construct Salem sets of arbitrary dimension $0<\alpha<1$;
the necessary modification would involve letting 
$N_j$ and $t_j$ be slowly increasing sequences such that
$\frac{\log t_j}{\log N_j}\to\alpha$.  We choose not to do
so here, as it would complicate the argument without contributing new ideas.

The constant $K$ will be set to equal $2^N$, but any other sufficiently
rapidly increasing function of $N$ would do.  The only purpose of
this constant is to ensure the uniformity of the constants $C_1,C_2$
in (A)-(B).  A reader interested only in a construction of Salem sets
in the traditional sense may as well set $K=1$.

The construction proceeds by iteration.  
Let $N_0$ and $t_0$ be integers such that $1\leq t_0\leq N_0$, and let
$\alpha= \frac{\log t_0}{\log N_0}$. 
Let $N=N_0^n$ and $t=t_0^n$, where $n$ is a large integer to be 
chosen later.
This will allow us to choose $N$ sufficiently
large depending on $\alpha$, while maintaining the relation
$\alpha=\frac{\log t}{\log N}$.  Let also
$K=2^N$, $M_j=KN^j$, $T_j=Kt^j$. 

We will construct inductively a sequence of sets $A_0,A_1,A_2,\dots$ such that 
$$A_0=\Big\{0,\frac{1}{K},\frac{2}{K},\dots,\frac{K-1}{K}\Big\},$$
$$A_{j+1}=\bigcup_{a\in A_j}A_{j+1,a},$$
where
\begin{equation}\label{ex-15a}
A_{j+1,a}\subset A^*_{j+1,a}:=\Big\{a,a+\frac{1}{M_{j+1}},a+\frac{2}{M_{j+1}},\dots,
a+\frac{N-1}{M_{j+1}}\Big\},
\end{equation}
and $\# A_{j+1,a} =t$ for each $a\in A_j$. In particular,
$\#A_j= T_j$.
Given such $A_j$, we define 
$$
E_j=\bigcup_{a\in A_j}[a,a+M_j^{-1}],
$$
$$E=\bigcap_{j=1}^\infty E_j.$$
Clearly, $E_1\supset E_2\supset E_3\supset\dots$, hence $E$ is a closed non-empty set.

There is a natural measure $\sigma$ on $E$, defined as follows.
Let $\calb_j$ be the family of all intervals of the form $I=[a,a+M_j^{-1}]$, $a\in A_j$,
and let $\calb=\bigcup \calb_j$.
For $F\subset E$, let
\begin{equation}\label{ex-11}
\sigma(F)=\inf\Big\{\sum_{i=1}^\infty T^{-1}_{j(i)}:\ F\subset\bigcup_{i=1}^\infty I_i,\ I_i\in\calb_{j(i)}\Big\}.
\end{equation}
Then $\sigma$ is the weak limit of the absolutely continuous measures
$\sigma_j$ with densities  
$$\phi_j=\sum_{a\in A_j}T_j^{-1}M_j{\bf 1}_{[a,a+M_j^{-1}]}.$$
In particular, we have
\begin{equation}\label{ex-12}
\sigma(I)=T^{-1}_{j}\hbox{ for all }I\in\calb_j.
\end{equation}

\begin{lemma}
Given $0<\alpha<1$ and $C_1 >1$, $\sigma$ satisfies the assumption (A) if 
$N=N(\alpha,C_1)$ has been chosen large enough. 
\end{lemma}
\begin{proof}
It suffices to prove that (A) holds if the interval $J=(x,x+\epsilon)$ is contained in
$[0,1]$.  Let $m_0$ be a positive integer such that $m_0^{-1}(m_0+1)<C_1$.
If $|J|\geq m_0 K^{-1}$, then there is an integer $m\geq m_0$
such that $mK^{-1}\leq |J|\leq (m+1)K^{-1}$.  We then have
$$
\sigma(J)\leq (m+1)K^{-1}\leq(m+1)\frac{|J|}{m}\leq C_1|J|\leq C_1|J|^\alpha,
$$
as required. If on the other hand $|J|\leq m_0K^{-1}$, 
let $m = m(J)$ be an integer such that $M_{m+1}^{-1} \leq |J| < m_0 M_m^{-1}$. 
Then $J$ is covered by at most $m_0+1$ intervals in $\mathcal B_m$, so (\ref{ex-12}) yields 
\[ \sigma(J \cap E)  \leq (m_0+1)T_m^{-1}. \]
Condition (A) will follow if we verify that $(m_0+1)T_m^{-1} \leq C_1 M_{m+1}^{-\alpha}$,
i.e.
$$(m_0+1)2^{-N} N^{-j\alpha}\leq C_1 2^{-N\alpha}N^{-(j+1)\alpha}.$$
But this simplifies to $2^{N(\alpha-1)}N^{\alpha}\leq C_1(m_0+1)^{-1}$,
which holds for any
$\alpha<1$ and $C_1>1$ if $N$ has been chosen large enough.  
\end{proof}

We must now prove that the sets $A_j$ can be chosen so that $\sigma$ 
obeys (B).  The proof will rely on probabilistic arguments inspired by those of Green in
\cite[Lemmas 14 and 15]{green-sumsets}; in particular, Lemma \ref{ex-lemma1} 
is almost identical to Lemma 14 in \cite{green-sumsets}.

For a finite set $A\subset \rr$, we will write
$$S_A(k)=\sum_{a\in A}e^{-2\pi ika}.$$

\begin{lemma}\label{ex-lemma1}
Let $B^*=\{0,\frac{1}{MN},\frac{2}{MN},\dots,\frac{N-1}{MN}\}$, where $M,N$ are large
integers, and let $t\in \{1,\dots,N\}$. Let
\begin{equation}\label{ex-2}
\eta^2t  = 32\ln(8MN^2).
\end{equation}
Then there is a set $B\subset B^*$ with $\#B=t$ such that
\begin{equation}\label{ex-3}
\Big|\frac{S_{B_x}(k)}{t}-\frac{S_{B^*}(k)}{N}\Big| \leq\eta\hbox{ for all }
k\in\zz,\ x=0,1,\dots,N-1,
\end{equation}
where
$$B_x=\Big\{\frac{(x+y)(\hbox{mod}\,N)}{MN}:\ \frac{y}{MN}\in B\Big\}.$$
\end{lemma}

The proof will be based on the following version of Bernstein's inequality,
which we also borrow from \cite{green-sumsets}.  We state it here
for completeness.  

\begin{lemma}\label{ex-lemma2}
Let $X_1,\dots,X_n$ be independent random variables with $|X_j|\leq 1$,
$\ee X_i=0$ and $\ee |X_j|^2=\sigma_j^2$.  Let $\sum \sigma_j^2\leq\sigma^2$,
and assume that $\sigma^2\geq 6n\lambda$.  Then 
\begin{equation}\label{ex-6}
\pp\Big(\Big|\sum_1^n X_j\Big|\geq n\lambda \Big)\leq 4e^{-n^2\lambda^2/8\sigma^2}.
\end{equation}
\end{lemma}

\noindent
{\em Proof of Lemma \ref{ex-lemma1}.} 
Let $B\subset B^*$ be a random set created by
choosing each $b\in B^*$ independently with probability $p=t/N$, where $t\in\nn$.
For each $b\in B^*$,
define the random variable $X_b(k)=(B(b)-p)e^{-2\pi ibk}$, where
we use $B(\cdot)$ to denote the characteristic function of the set $B$.
Then $X_b(k)$ obey the assumptions of Lemma \ref{ex-lemma2} with 
$\sigma_b^2=\ee|X_b(k)|^2=\ee|B(b)-p|^2=p-p^2\in(p/2,p)$, so that 
$t/2\leq \sigma^2\leq t$.  We also have 
$$
\frac{S_{B}(k)}{t}-\frac{S_{B^*}(k)}{N}=t^{-1}\sum X_b(k)$$
and in particular $\sum X_b(0)=\#B -t.$
Applying Lemma \ref{ex-lemma2} with
$n=N$ and $\lambda=\eta p/2$, we see that 
\begin{equation}\label{ep-1}
\pp\Big(\Big|
\frac{S_{B}(k)}{t}-\frac{S_{B^*}(k)}{N}\Big|  \geq\eta/2 \Big)
\leq 4\exp(-\eta^2t/32).
\end{equation}
By the same argument, (\ref{ep-1}) holds with $B$ replaced by $B_x$.

Note that
$S_{B^*}(k)$ and $S_{B_x}(k)$ are periodic with period $MN$, hence it suffices 
to consider $k\in\{0,1,\dots,MN-1\}$.  
Thus the probability that the above event occurs for any such $k$
and with $B$ replaced by any $B_x$ is bounded by $4MN^2\exp(-\eta^2 t/32)$, 
which is less than $1/2$ whenever (\ref{ex-2}) holds. Thus with probability at least
$1/2$ we have
$$\Big|\frac{S_{B_x}(k)}{t}-\frac{S_{B^*}(k)}{N}\Big| \leq\frac{\eta}{2}\hbox{ for all }
k\in\zz,\ x=0,1,\dots,N-1.$$
Note further that the last inequality
with $k=0,x=0$ implies that $|\#B-t|\leq \eta t/2$.  Modifying $B$ by
at most $\eta t/2$ elements, we get a set of cardinality $t$ obeying (\ref{ex-3}).
This proves the lemma.
\hfill$\square$

\bigskip

\begin{lemma}\label{ex-lemma3}
The sets $A_j$ can be chosen so that
\begin{equation}\label{ex-5}
\Big|\widehat{\phi_{j+1}}(k)-\widehat{\phi_j}(k)\Big|\leq 
16\min\Big(1,\frac{M_{j+1}}{|k|}\Big)T_{j+1}^{-1/2}\ln(8M_{j+1}).
\end{equation}
\end{lemma}

\begin{proof}
Since the index $j$ will be fixed throughout this proof, we drop it from the notation and write
$A=A_j$, $A'=A_{j+1}$, $T=T_j$, $M=M_j$.  
With this notation,
we have
$\phi_j=\sum_{a\in A}T^{-1}M{\bf 1}_{[a,a+M^{-1}]}$, 
hence
\begin{equation}\label{ex-1}
\widehat{\phi_j}(k)=MT^{-1}\sum_{a\in A} \int_a^{a+M^{-1}}e^{-2\pi ikx}dx
=\frac{1-e^{2\pi ik/M}}{2\pi ik/M}T^{-1}S_A(k). 
\end{equation}
Let also $B^*,B,B_x$ be as in Lemma \ref{ex-lemma1}, and let
$$A'=\bigcup_{a\in A}(a+B_{x(a)}),$$
where $x(a)$ is chosen randomly from the set $\{0,1,\dots,N-1\}$ and takes each value 
with probability $N^{-1}$, and the choices are independent for different $a$'s.
Then
\begin{equation}\label{ex-101}
\widehat{\phi_{j+1}}(k)-\widehat{\phi_j}(k)
=\frac{1-e^{2\pi ik/MN}}{2\pi ik/MN}\cdot T^{-1}\sum_{a\in A}\calx_a(k),
\end{equation}
where 
$$
\calx_a(k)=\frac{S_{B_{x(a)}+a}(k)}{t}-\frac{S_{B^*+a}(k)}{N}.
$$

Fix $k$.
We consider $\calx_a(k)$, $a\in A,$ as independent random variables.  It is easy to check that
$\ee \calx_a(k)=0$ for each $a$.  By Lemma \ref{ex-lemma1}, $|\calx_a(k)|\leq\eta$.
Applying Lemma \ref{ex-lemma2} to $\calx_a(k)$, with $n=T$ and $\sigma^2=T\eta^2$,
we find that
\begin{equation}\label{ep-2}
\pp\Big(\Big|
T^{-1}\sum\calx_a(k)\Big|  \geq\lambda \Big)
\leq 4\exp(-\lambda^2T/8\eta^2).
\end{equation}
Thus the probability that this happens for any $k\in\{0,1,\dots,MN-1\}$ is bounded by
$4MN\exp(-\lambda^2T/8\eta^2)$, which is less than $1/2$ if
$\lambda^2 T\geq 4\eta^2\ln(8MN)$.  If $\eta$ is as in (\ref{ex-2}), it is easy to check that
the last inequality holds for
$\lambda=16 (Tt)^{-1/2}\ln(8MN)$.  This together with (\ref{ex-101}) completes the proof.

\end{proof}

\begin{lemma}\label{ex-lemma10}
Let $C_2>0$ and $0<\beta<\alpha<1$.  If $N$ is large enough, depending
on $C_2,\alpha,\beta$, the sets $A_j$ can be chosen so that (B) holds.
\end{lemma}

\begin{proof}
Let $A_j$ be as in Lemma \ref{ex-lemma3}.
By (\ref{ex-5}), it suffices to prove that for $N$ sufficiently large we have
\begin{equation}\label{ex-50}
\sum_{j=1}^\infty 
\min\Big(1,\frac{M_{j}}{|k|}\Big)T_{j}^{-1/2}\ln(8M_{j})<\frac{C_2}{16} |k|^{-\beta/2},\ k\neq 0.
\end{equation}
We may assume that $k>0$.  Plugging in the values of $M_j$ and $T_j$, we see that
(\ref{ex-50}) is equivalent to
\begin{equation}\label{ex-51}
\sum_{j=1}^\infty 
\min\Big(1,\frac{2^N N^j}{k}\Big)2^{-N/2}N^{-j\alpha/2}
(N\ln 2+\ln 8+j\ln N)<\frac{C_2}{16} k^{-\beta/2}.
\end{equation}
We write the sum in (\ref{ex-51}) as
$$
\sum_{j=1}^\infty 
\min\Big(1,\frac{2^N N^j}{k}\Big)2^{-N/2}N^{-j\beta/2}N^{-j(\alpha-\beta)2}
(N\ln 2+\ln 8+j\ln N).
$$
It is a simple exercise in calculus to check that $N^{-j(\alpha-\beta)/2}j\ln N\leq
2(\alpha-\beta)^{-1}$.
Hence the sum in (\ref{ex-51}) is bounded by 
$$
\sum_{j=1}^\infty 
\min\Big(1,\frac{2^N N^j}{k}\Big)2^{-N/2}N^{-j\beta/2}
(N\ln 2+\ln 8+2(\alpha-\beta)^{-1}).
$$
We write the last sum as $S_1+S_2$, where $S_1$ is the sum over all
$j$ with $1\leq j\leq \frac{\ln k-N\ln 2}{\ln N}$ and $S_2$ is the sum over all remaining
values of $j$.  We first estimate $S_1$.  We have
$$
S_1=2^{N/2}k^{-1}(N\ln 2+\ln 8+2(\alpha-\beta)^{-1})
\sum_{1\leq j\leq \frac{\ln k-N\ln 2}{\ln N}}N^{j(1-\frac{\beta}{2})}.
$$
The last sum is bounded by
$$
2N^{\frac{\ln k-N\ln 2}{\ln N}(1-\frac{\beta}{2})}
=2k^{1-\frac{\beta}{2}}2^{-(1-\frac{\beta}{2})N},
$$
hence
\begin{align*}
S_1
&\leq 2^{N/2}k^{-1}(N\ln 2+\ln 8+2(\alpha-\beta)^{-1})
\cdot 2k^{1-\frac{\beta}{2}}2^{-(1-\frac{\beta}{2})N}\\
&=2^{\frac{\beta-1}{2}N+1}(N\ln 2+\ln 8+2(\alpha-\beta)^{-1})
k^{-\beta/2},
\end{align*}
which is bounded by $C_2k^{-\beta/2}/100$ if $N$ is large enough, depending on $C_2$,
$\alpha,\beta$. 

We now turn to $S_2$:
$$
S_2=2^{-N/2}(N\ln 2+\ln 8+2(\alpha-\beta)^{-1})
\sum_{j > \frac{\ln k-N\ln 2}{\ln N}}N^{-j\beta/2},
$$
and the last sum is bounded by
$$
N^{-\frac{\ln k-N\ln 2}{\ln N}\frac{\beta}{2}}
=k^{-\beta/2}2^{-\beta N/2}.
$$
Thus
\begin{align*}
S_2
&\leq 2^{-N/2}(N\ln 2+\ln 8+2(\alpha-\beta)^{-1})
\cdot k^{-\frac{\beta}{2}}2^{-\beta N/2}\\
&=2^{\frac{\beta-1}{2} N}(N\ln 2+\ln 8+2(\alpha-\beta)^{-1})
k^{-\beta/2},
\end{align*}
which  again is bounded by $C_2k^{-\beta/2}/100$ if $N$ is large enough.  
This ends the proof of the lemma.

\end{proof}

{\bf Remark.} The same argument shows that if $N$ and $t$ are fixed,
then $\sigma$ obeys (B)
for all $\beta<\alpha$ with some constant $C_2=C_2(\alpha, \beta)$. 
Thus the sets constructed here are also Salem sets in the traditional sense.


\section{Salem's Construction}
The purpose of this section is to establish that Salem's random
construction of Salem sets \cite{salem} provides a rich class of
examples for which conditions (A) and (B) can be verified. More precisely,
we will show that there exist absolute large constants $C_1$ and $C_2$
such that for $\alpha$ arbitrarily close to 1, there are sets occurring
with high probability that satisfy condition (A) with exponent $\alpha$
and constant $C_1$, and (B) with some $\beta > 2/3$, $C_2 > 0$
and $B = \frac{1}{2}$. It is important for our analysis that the
constants $C_1$, $C_2$ remain bounded as $\alpha \rightarrow 1$.   

Let us recall Salem's construction of these sets, which is based
on a generalization of the Cantor construction. Given an integer
$d \geq 2$, let $0 < a_1 < a_2 < \cdots a_d < 1$ be $d$ numbers that
are linearly independent over the rationals. Let $\kappa > 0$ be a
number satisfying      
\begin{equation} 0 < \kappa < \min \{a_j - a_{j-1} : 1 \leq j \leq d\} \quad \text{ and } \quad \kappa < 1 - a_d. \label{kappa-condition}  \end{equation} 
Given an interval $[a,b]$ of length $L$, a dissection of type $(d, a_1, a_2, \cdots, a_d, \kappa)$ is performed on $[a,b]$ by calling each of the closed intervals $[a + La_j, a+L(a_j + \kappa)], 1 \leq j \leq d$, white and the complementary intervals black. 

Let us fix the numbers $d$, $a_1, \cdots, a_d$, and an infinite
sequence $\{ \kappa_m : m \geq 1\}$, each of whose elements satisfies
(\ref{kappa-condition}). Starting with $E_0 = [0,1]$, we perform a
dissection of type $(d, a_1, \cdots, a_d, \kappa_1)$ and remove the
black intervals, thereby obtaining a set $E_1$ which is a union of $d$ intervals each of length $\kappa_1$. On each of the component intervals of $E_1$, we perform a dissection of type $(d, a_1, \cdots, a_d, \kappa_2)$, remove the black intervals and so obtain a set $E_2$ of $d^2$ intervals each of length $\kappa_1\kappa_2$. After $n$ steps we obtain a set $E_n$ of $d^n$ intervals, each of length $\kappa_1 \cdots \kappa_n$. Letting $n \rightarrow \infty$, we obtain a perfect nowhere dense set $E = \cap_{n=1}^{\infty} E_n$, which has Lebesgue measure zero if $d^n \kappa_1 \cdots \kappa_n \rightarrow 0$. 

For each $n \in \mathbb N$, let $F_n$ be a continuous nondecreasing function satisfying 
\begin{itemize}
\item $F_n(x) = 0$ for $x \leq 0$; $F_n(x) = 1$ for $x \geq 1$. 
\item $F_n$ increases linearly by $d^{-n}$ on each of the $d^n$ white intervals constituting $E_n$.
\item $F_n$ is constant on every black interval complementary to $E_n$.   
\end{itemize}    
The pointwise limit $F = \lim_{n \rightarrow \infty} F_n$ is a nondecreasing continuous function with $F(0) = 0$, $F(1) = 1$, and can therefore be realized as the distribution function of a probability measure $\mu$. The Fourier transform of $\mu$ is given by 
\begin{equation} \widehat{\mu}(\xi) = P(\xi) \prod_{n=1}^{\infty} P(\xi \kappa_1 \cdots \kappa_n),
\quad \text{ where } \quad P(\xi) =  \frac{1}{d} \sum_{j=1}^{d} e^{2 \pi i a_j \xi}. \label{fourier-coeff} \end{equation}

Given $\alpha < 1$ (note that $d$ and $\alpha$ are independent parameters), let us set $\kappa = d^{-\frac{1}{\alpha}}$, and further restrict our choice of $a_j$ so that they satisfy 
\begin{equation} 0 <  a_1 < \frac{1}{d} - \kappa, \quad \text{ and } \quad \kappa < a_j - a_{j-1} < \frac{1}{d} \text{ for } 2 \leq j \leq d. \label{revised-a}  \end{equation} 
Let $\Xi = \Xi(d, \alpha)$ be the collection of all infinite sequences $\mathbf k = \{\kappa_m : m \geq 1 \}$ satisfying
\[ \Bigl(1 - \frac{1}{2m^2}\Bigr)\kappa \leq \kappa_m \leq \kappa, \quad m \geq 1. \]  Repeating the construction outlined in the previous paragraph with a fixed choice of $(a_1, \cdots, a_d)$ as in (\ref{revised-a}) and different choices of $\mathbf k \in \Xi$, we obtain an uncountable collection of sets $E = E[\mathbf k]$, all of which have Hausdorff dimension $\alpha$. In fact, the supporting measures $\mu = \mu[\mathbf k]$ satisfy a ball condition of the form (A), a fact that was observed in \cite{mock}. The following result is a rephrasing of Proposition 3.2 in \cite{mock}, with special attention to the implicit constants.  
\begin{proposition}[\cite{mock}]
There exists an absolute constant $C_0$ such that for every $\alpha < 1$ and $\mathbf k \in \Xi$, the corresponding measure $\mu = \mu[\mathbf k]$ satisfies 
\[ \mu[x, x+r] \leq dC_0 r^{\alpha}, \qquad 0< r \leq 1. \] 
Thus condition (A) holds with $C_1 = dC_0$.  
\end{proposition}

We now turn to (B). One of the main results in \cite{salem} is that for fixed $d \geq 2$, $\alpha < 1$ and $(a_1, \cdots, a_d)$ satisfying (\ref{revised-a}), there exists a parametrization of $\Xi = \{\mathbf k(t) : 0 \leq t \leq 1 \}$ such that for almost every $t \in [0,1]$, $\widehat{\mu}[\mathbf k(t)]$ satisfies a decay condition. We need a stronger version of this result that formalizes how the implicit constant in the Fourier decay condition depends on $\alpha$. The main result in this section is the next proposition.  
\begin{proposition}
For all $d$ sufficiently large, there exist constants $C_2 \geq 1$, $\epsilon_0 \ll 1$ and $\beta > \frac{4}{5}$ (depending on $d$) with the following property. For every $\alpha \in (1 - \epsilon_0, 1)$,  there exist numbers $(a_1, \cdots, a_d)$ satisfying (\ref{revised-a}) such that the random measures $\mu = \mu[\mathbf k]$ constructed by Salem based on $(a_1, \cdots, a_d)$ obey the Fourier decay estimate  
\[ |\widehat{\mu}(\xi)| \leq \frac{C_2}{(1 - \alpha)^{\frac{1}{2}}} |\xi|^{- \frac{\beta}{2}} \quad \text{ for all } \xi \ne 0,   \] 
with large probability.   \label{salem-mainprop}
\end{proposition}
The proof of the proposition is based on the following three lemmas. 
\begin{lemma}
Given any $m, M\geq 10$, there exists $\mathbf x = (x_1, \cdots, x_m) \in (0,1)^m$ such that 
\[ |\mathbf x \cdot \mathbf r| \geq M^{-2m} \quad \text{ for all } \; 0 \ne \mathbf r \in \mathbb Z^m, \; ||\mathbf r||_{\infty} \leq M.  \] \label{ingredient-lemma1}
\end{lemma}
\begin{proof}
The proof is a simple volume estimation argument. Given any $\mathbf r$ as in the statement of the lemma, let 
\[ V_{\mathbf r} = \Bigl\{\mathbf x \in (0,1)^m : | \mathbf x \cdot \mathbf r | < \epsilon ||\mathbf r||_{\infty} \Bigr\}. \] 
Then $|V_{\mathbf r}| \leq \epsilon$. Since the number of possible choices of $\mathbf r$ is $(2M+1)^m - 1$, $|\cup \{V_{\mathbf r} : 0 \ne \mathbf r \in \mathbb Z^m, ||\mathbf r||_{\infty} \leq M  \}| \leq (2M+1)^m \epsilon$. Choosing any $\epsilon < (2M+1)^{-m}$ would therefore guarantee the existence of $\mathbf x\in [0,1]^m \setminus \cup V_{\mathbf r}$. In particular, $\epsilon = M^{-2m}$ suffices. 
\end{proof}
\begin{lemma}
Let $Q(\xi) = \sum_{j=1}^d \lambda_j e^{2 \pi i b_j \xi}$, where $\mathbf b = (b_1, \cdots, b_d)$ is any collection of real numbers linearly independent over the rationals. Given $s > 0$, let 
\begin{equation} \delta_s(\mathbf b) = \inf \Bigl\{|\mathbf b \cdot \mathbf j| : 0 \ne \mathbf j \in \mathbb Z^d,\; \mathbf j \cdot \mathbf 1  = 0, \; ||\mathbf j||_{\infty} \leq \frac{s}{2}+1  \Bigr\} > 0, \label{delta-def} \end{equation} where $\mathbf 1 = (1, \cdots, 1) \in \mathbb R^d$. Then there exists a constant $c = c(d,s) > 0$ such that for all $T \geq T_0 = c(d,s)(\delta_s(\mathbf b))^{-1}$ and all $t \in \mathbb R$:
\[ \frac{1}{T} \int_{t}^{T+t} |Q(\xi)|^s d\xi \leq 2 \left(\frac{s}{2} + 1 \right)^{\frac{s}{2}} \bigl(\sum_{j=1}^{d} \lambda_j^2 \bigr)^{\frac{s}{2}}. \] \label{salem-lemma}
\end{lemma}
\begin{proof}
This result is a variant of the lemma in section 3 of \cite{salem}. Unlike \cite{salem} however, 
we are mainly concerned with the explicit dependence of $T_0$ on $\delta_s(\mathbf b)$, so we revise the proof with attention to this detail. Let $2q$ be the even integer such that $s \leq 2q < s+2$. Then $|Q(\xi)|^{2q} = Q(\xi)^q \overline{Q}(\xi)^q = Q_1 + Q_2(\xi)$, where 
\begin{align*}
Q_1 &= \sum_{\begin{subarray}{c}\mathbf j \geq 0 \\ \mathbf j \cdot \mathbf 1 = q \end{subarray}} |\pmb{\lambda}^{2 \mathbf j}| \left(\frac{q!}{\mathbf j!} \right)^2 < q! \sum_{\begin{subarray}{c} \mathbf j \geq 0 \\ \mathbf j \cdot \mathbf 1  = q \end{subarray}} |\pmb{\lambda}^{2 \mathbf j}| \frac{q!}{\mathbf j!} = q! \Bigl(\sum_{j=1}^d |\lambda_j|^2 \Bigr)^{q}, \text{ and } \\ Q_2(\xi) &= \sum_{\begin{subarray}{c} \mathbf j \ne \mathbf i \\ \mathbf j \cdot \mathbf 1  = \mathbf i \cdot \mathbf 1  = q  \end{subarray}} \pmb{\lambda}^{\mathbf j} \overline{\pmb{\lambda}}^{\mathbf i} e^{2 \pi i \xi (\mathbf j - \mathbf i)\cdot \mathbf b}. 
\end{align*} 
Here $\pmb{\lambda} = (\lambda_1, \cdots, \lambda_d)$; the vectors $\mathbf j = (j_1, \cdots, j_d)$ and $\mathbf i = (i_1, \cdots, i_d)$ are multi-indices consisting of non-negative integer entries; $\mathbf j! = j_1! \cdots j_d!$, $\pmb{\lambda}^{\mathbf j} = \lambda_1^{j_1} \cdots \lambda_d^{j_d}$. Thus
{\allowdisplaybreaks \begin{align*}
\frac{1}{T} \int_{t}^{T+t} |Q(\xi)|^{2q} \, d\xi &\leq q! \Bigl(\sum_{j=1}^{d} |\lambda_j|^2 \Bigr)^q + \frac{1}{T} \left|\int_{t}^{T+t} Q_2(\xi) \, d\xi \right| \\ &\leq q! \Bigl(\sum_{j=1}^{d} |\lambda_j|^2 \Bigr)^q + \Bigl| \sum_{\begin{subarray}{c} \mathbf j \ne \mathbf i \\ \mathbf j \cdot \mathbf 1  = \mathbf i \cdot \mathbf 1 = q  \end{subarray}} \pmb{\lambda}^{\mathbf j} \overline{\pmb{\lambda}}^{\mathbf k} \frac{\left[e^{2 \pi i (\mathbf j - \mathbf i) \cdot \mathbf b (T+t)} - e^{2 \pi i (\mathbf j - \mathbf i) \cdot \mathbf b t} \right]}{T(\mathbf j - \mathbf i) \cdot \mathbf b} \Bigr|  \\ &\leq q! \Bigl(\sum_{j=1}^{d} |\lambda_j|^2 \Bigr)^q + c(d,q) ||\pmb{\lambda}||_{\infty}^{2q} \frac{1}{\delta_s(\mathbf b)T}, 
\end{align*}
and hence $\leq 2q! \left(\sum_{j=1}^{d} |\lambda_j|^2 \right)^{q}$
if $T \geq T_0 = c(d, q)(\delta_s(\mathbf b))^{-1}$. By H\"older's inequality, 
\begin{align*}
\left[ \frac{1}{T} \int_{t}^{T+t} |Q(\xi)|^s \, d\xi  \right]^{\frac{1}{s}} &\leq \left[ \frac{1}{T} \int_{t}^{T+t} |Q(\xi)|^{2q} \, d\xi  \right]^{\frac{1}{2q}} \\ &\leq (2q!)^{\frac{1}{2q}} \Bigl(\sum_{j=1}^{d}|\lambda_j|^2 \Bigr)^{\frac{1}{2}} \\ &\leq 2^{\frac{1}{2q}} q^{\frac{1}{2}} \Bigl( \sum_{j=1}^{d} |\lambda_j|^2 \Bigr)^{\frac{1}{2}} \\ &\leq 2^{\frac{1}{s}} \Bigl(\frac{s}{2} + 1 \Bigr)^{\frac{1}{2}} \Bigl( \sum_{j=1}^{d} |\lambda_j|^2 \Bigr)^{\frac{1}{2}}, 
\end{align*}}
whence the result follows. 
\end{proof}
\begin{lemma}
Given $d$ and $s \geq 2$, there exist positive constants $c_0$ and $\epsilon_0 \ll 1$ depending only on these parameters such that for all $\alpha \in (1 - \epsilon_0, 1)$, there is $\mathbf a = (a_1, \cdots, a_d)$ satisfying (\ref{revised-a}) and $\delta_s(\mathbf a) \geq c_0(1 - \alpha)$. Here $\delta_s(\cdot)$ is as in (\ref{delta-def}). \label{ingredient-lemma3} 
\end{lemma}
\begin{proof}
We replace the variables $\mathbf a = (a_1, \cdots, a_d)$ satisfying (\ref{revised-a}) by the new set $\pmb{\eta} = (\eta_2, \cdots, \eta_d) \in [0,1]^{d-1}$, with $\pmb{\zeta} = (\zeta_2, \cdots, \zeta_d)$ being an intermediate set of coordinates. These are defined as follows:
\begin{align*}
\zeta_j &= d^{\frac{1}{\alpha}} (a_j - a_{j-1}) \; \text{ for } 2 \leq j \leq d, \text{ so that } 1 < \zeta_j < d^{\frac{1}{\alpha} - 1}, \quad \text{ and } \\
\eta_j &= \frac{\zeta_j -1}{d^{\frac{1}{\alpha}-1} - 1} \; \text{ for } 2 \leq j \leq d,  \text{ so that } 0 < \eta_j < 1. 
\end{align*}  
For any $0 \ne \mathbf j \in \mathbb Z^d$ with $\mathbf j \cdot \mathbf 1 = 0$ and $||\mathbf j||_{\infty} \leq \frac{s}{2} + 1$, the linear functional $\mathbf a \mapsto \mathbf j \cdot \mathbf a$ may be expressed in these new coordinates as: 
\begin{align*}
\mathbf j \cdot \mathbf a &= \left(\frac{1}{d}\right)^{\frac{1}{\alpha}} \left[\zeta_2(j_2 + \cdots + j_d) + \zeta_3(j_3 + \cdots + j_d) + \zeta_d j_d \right] \\ &= \left(\frac{1}{d}\right)^{\frac{1}{\alpha}} [\zeta_2 m_2 + \cdots + \zeta_d m_d] \\ &= \left(\frac{1}{d}\right)^{\frac{1}{\alpha}} \left[\mathbf m \cdot \mathbf 1 + \left(d^{\frac{1}{\alpha} - 1} - 1 \right) \pmb{\eta} \cdot \mathbf m\right],
\end{align*} 
where $\mathbf m = (m_2, \cdots, m_d) \in \mathbb Z^{d-1}$, $m_{\ell} = j_{\ell} + j_{\ell + 1} + \cdots + j_d$, so that $0 < ||\mathbf m||_{\infty} \leq ds$. Thus to prove the lemma it suffices to show that for all $\alpha$ sufficiently close to 1, there exists $\pmb{\eta} \in (0,1)^{d-1}$ satisfying
\begin{equation} 
\begin{aligned} 
\inf \Bigl\{ \Bigl(\frac{1}{d}\Bigr)^{\frac{1}{\alpha}} \left| \mathbf m \cdot \mathbf 1 + \left( d^{\frac{1}{\alpha} - 1} - 1 \right) \pmb{\eta} \cdot \mathbf m \right| : 0 \ne \mathbf m \in \mathbb Z^{d-1}, ||\mathbf m||_{\infty} \leq ds \Bigr\} \\ \geq c_0 (1 - \alpha). 
\end{aligned}  \label{toshow} \end{equation}
We consider two cases. If $\mathbf m \cdot \mathbf 1 \ne 0$, then for every $\pmb{\eta} \in (0,1)^{d-1}$, 
\[\Bigl(\frac{1}{d}\Bigr)^{\frac{1}{\alpha}} \left| \mathbf m \cdot \mathbf 1 + \left( d^{\frac{1}{\alpha} - 1} - 1 \right) \pmb{\eta} \cdot \mathbf m \right| \geq \Bigl(\frac{1}{d}\Bigr)^{\frac{1}{\alpha}} \left[ 1 - \left(d^{\frac{1}{\alpha}} - 1 \right) d^2 s \right] \geq \frac{1}{2d^2}, \]
provided $\alpha$ is close enough to 1 to ensure $(d^{\frac{1}{\alpha} - 1} - 1) d^2 s \leq  2(1 - \alpha) d^3 s \log d \leq \frac{1}{2}$. This of course gives a better estimate than required by (\ref{toshow}). If $\mathbf m \cdot \mathbf 1 = 0$, then by Lemma \ref{ingredient-lemma1} there exists $\pmb{\eta} \in (0,1)^{d-1}$ such that 
\[ \Bigl(\frac{1}{d}\Bigr)^{\frac{1}{\alpha}} \bigl(d^{\frac{1}{\alpha} - 1} - 1\bigr) |\pmb{\eta} \cdot \mathbf m| \geq d^{-2} (1 - \alpha) \log d |\pmb{\eta} \cdot \mathbf m| \geq d^{-2} (1 - \alpha) \log d (ds)^{-2(d-1)}.\]
This completes the proof of (\ref{toshow}) and hence the lemma.  
\end{proof}
\begin{proof}[\bf{Proof of Proposition \ref{salem-mainprop}}]
We follow the proof of Theorem II in \cite{salem} with minor modifications
and special attention to constants. Let $\epsilon > 0$ be sufficiently small
(for instance any $\epsilon < \frac{1}{10}$ suffices). The value of
$\epsilon$ will remain fixed throughout the proof and the 
constant $C_2$ in condition (B) will depend on this choice of $\epsilon$.
Let $s = 3/\epsilon$, and $d$ be the smallest integer $\geq 2$ such that 
\[ \sqrt{d} \geq 2 \left(\frac{s}{2} + 1 \right)^{\frac{s}{2}}.  \]
By Lemma \ref{ingredient-lemma3}, there exist constants $c_0$ and
$\epsilon_0 > 0$ depending only on $\epsilon$ such for every $\alpha
\in (1 - \epsilon_0, 1)$ we can find $\mathbf a = (a_1, \cdots, a_d)$
satisfying (\ref{revised-a}) and $\delta_s(\mathbf a) \geq c_0 (1 -
\alpha)$. For this choice of $\mathbf a$ and $s$, and $P$ as in (\ref{fourier-coeff}), the conclusion of Lemma  \ref{salem-lemma} holds with $Q$ replaced by $P$, and $T_0 = C_{\epsilon}/(1 - \alpha)$ for some large constant $C_{\epsilon}$ (independent of $\alpha$).  Let $\mu_t$ denote the random measure generated by the sequence $\mathbf k(t) \in \Xi$, $0 \leq t \leq 1$, and the fixed choice of $\mathbf a$ above. Choosing $\theta = 2(\frac{s}{2}-1)/(s-1)$, the proof of Theorem II in \cite{salem} yields 
\[ \int_{0}^{1} |\widehat{\mu}_t(n)|^s dt \leq \frac{1}{|n|^{\alpha \left(\frac{s}{2} - 1 \right)}}, 
\quad \text{ for } |n| \geq T_0.   \]
Writing $\alpha(\frac{s}{2} -1) = 2 +\gamma$, we obtain for some absolute constant $C > 0$, 
\begin{equation} \sum_{n \geq T_0} n^\gamma \int_0^1 |\widehat{\mu}_t(n)|^s dt \leq \sum_{n \geq T_0}\frac{1}{n^2}   \leq \frac{C}{T_0}. \label{expectation}  \end{equation}
By (\ref{expectation}) and Chebyshev's inequality, 
\begin{align*}
\mathbb P \Bigl(\Bigl\{t : \sup_{n \ne 0} \bigl| |n|^{\frac{\gamma}{s}} 
\widehat{\mu}_t(n) \bigr| > N \Bigr\} \Bigr) 
&\leq \mathbb P \Bigl(\sum_{n \ne 0} |n|^{\gamma} \bigl|\widehat{\mu}_t(n) \bigr|^s > N^s \Bigr) \\ 
&\leq \mathbb P\Bigl( \sum_{n \geq T_0} |n|^{\gamma} |\widehat{\mu}_t(n)|^s > N^s - T_0^{\gamma+1}\Bigr) \\ &\leq \frac{C/T_0}{N^s - T_0^{\gamma+1}}.
\end{align*}     
Choosing $N = \sqrt{T_0}$, and observing that 
\[ \frac{\gamma}{s} = \frac{\alpha}{2} - \frac{(\alpha + 2)}{s} =
\frac{\alpha}{2} - \epsilon \frac{(\alpha + 2)}{3}, \]
we deduce that with large probability, condition (B) holds with some
exponent $\beta > \frac{4}{5}$, some constant $C_2$ (uniform in $\alpha$)
and $B = \frac{1}{2}$, provided that $\alpha$ is close to 1,
and $\epsilon$ is sufficiently small.  
\end{proof}

\section{Brownian images}
\label{sec-bm}

In this section, we address the question of finding 3-term arithmetic progressions in the
``random" Salem sets constructed by Kahane \cite{kahane}.  

Fix $\alpha\in(0,1)$, and let $F_{\alpha}$ be a subset of $[0,1]$ of Hausdorff dimension $\alpha/2$. 
Let $\theta = \theta_{\alpha}$ be a probability measure on a set $F_{\alpha}\subset\rr$ such that \begin{equation} \theta_{\alpha}(I) \leq C_0 |I|^{\frac{\alpha}{2}} \label{Frostman} \end{equation} 
for each interval $I$. By Frostman's lemma, such a measure exists provided $F_{\alpha}$ 
has positive Hausdorff measure of order $\frac{\alpha}{2}$. The constant $C_0$ can be chosen 
uniform in $\alpha$ as $\alpha \rightarrow 1$, provided the $\frac{\alpha}{2}$-dimensional 
Hausdorff measure of $F_{\alpha}$ remains bounded away from zero (see \cite{kahane}, p130). 
We will always assume this to be the case. Let $W(\cdot)$ denote the one-dimensional 
Brownian motion, and $\mu$ the image of $\theta$ by $W$, i.e.,  
\[ \int f \, d\mu  = \int_{0}^{1} f(W(t)) \, d\theta(t), \]
so that $\mu$ is a random measure on $\mathbb R$. The Fourier transform of $\mu$ is given by 
\[ \widehat{\mu}(\xi) = \int_{0}^{1} e^{-2 \pi i \xi W(t)} d\theta (t). \]
Kahane \cite{kahane} proves that under the above assumptions, 
$E=\supp\mu$ is almost surely a Salem set: for all $\beta<\alpha$ we have almost surely
\begin{equation}\label{br-e1}
\sup_{\xi\in\rr}|(1+|xi|)^{\frac{\beta}{2}} |\widehat{\mu}(\xi)| <\infty.
\end{equation}

It is not difficult to modify Kahane's argument so as to show that Assumption (B) holds (i.e.
the implicit constants in (\ref{br-e1}) can be chosen independent of $\alpha$) with probability
at least $1/2$.  On the other hand, Assumption (A) is not expected to hold: it is known in the 
probabilistic literature (see e.g. \cite{DPRZ}) that Brownian images of measures as in 
(\ref{Frostman}) obey an estimate 
similar to (A) but with an additional factor of $\log(\epsilon^{-1})$, and that this is optimal.
Thus our Theorem \ref{thm-main} does not apply in this case.  

Nonetheless, we are able to show that Brownian images contain non-trivial 3-term arithmetic
progressions with positive probability. Our approach here will not rely on Theorem 
\ref{thm-main}; instead, we will appeal directly to a variant of Proposition \ref{prop-measure}
which we now state.

\begin{proposition}\label{prop-measure2}
Let $\mu$ be a probability measure supported on a closed set $E\subset \rr$ 
such that (\ref{br-e1}) holds for some $\beta\in(2/3,1]$.
Assume furthermore that $\Lambda(\mu,\mu,\mu)>0$, where
\begin{equation}\label{br-e2}
\Lambda(\mu,\mu,\mu)=\int \whm^2(\xi)\whm(-2\xi)d\xi.
\end{equation}
Then there are $x,y\in E$ such that $x\neq y$ and $\frac{x+y}{2}\in E$.

\end{proposition}

This differs from the statement of Proposition \ref{prop-measure}
in that $\mu$ is no longer required to be compactly supported and that, accordingly,
the $\Lambda$ quantity is now defined using the continuous Fourier transform rather
than the Fourier series.  Proposition \ref{prop-measure2} follows by exactly the same
argument as in Section \ref{sec-measure}, except that $\nu$ will now be a linear functional
on $C_c(\rr^2)$, the space of all compactly supported continuous functions on $\rr^2$.

Our desired conclusion now follows from Proposition \ref{prop-measure2} and
the next proposition.

\begin{proposition}\label{prop-kahane2}
There is a constant $c>0$, depending only on the constant $C_0$ in 
(\ref{Frostman}) but independent of $\alpha$ as long as $\alpha>\frac{2}{3}+\epsilon$, such that
$$
\pp(\Lambda(\mu,\mu,\mu)>0)\geq c.
$$
\end{proposition}

The proof of Proposition \ref{prop-kahane2} will rely on the following
result due to Kahane \cite{kahane}, which we state here without proof.
(Strictly speaking, this is only proved in \cite{kahane} for integer $q$,
but the extension to all $q>0$ follows trivially from H\"older's
inequality.)

\begin{proposition}\label{prop-kahane}
\cite[pp. 254--255]{kahane}
Let $\mu$ be as defined above. Then for any $\xi \ne 0$ and any $q > 0$, 
\[ \mathbb E \left[ |\widehat{\mu}(\xi)|^{2q} \right] \leq (C'_0 q |\xi|^{-\alpha})^{q},  \] 
where $C'_0$ depends only on $C_0$ but not on $q$ or $\alpha$.
\label{kahane-result}
\end{proposition} 

\begin{lemma}
There is a constant $C$, depending only on $C_0$, such that for $\mu$ as above
we have
\begin{equation}\label{br-e10}
\mathbb E \left[ \int | \widehat{\mu}(\xi)|^2 |\widehat{\mu}(-2 \xi)| \, d\xi \right] <C,
\end{equation}
\begin{equation}\label{br-e11}
\mathbb E \left[\left( \int | \widehat{\mu}(\xi)|^2
|\widehat{\mu}(-2 \xi)| \, d\xi \right)^2\right] <C.
\end{equation}
In particular, this implies $\mathbb E(|\Lambda (\mu, \mu, \mu)|) < \infty$
and $\mathbb E(|\Lambda (\mu, \mu, \mu)|^2) < \infty$.
\label{exp-conv}
\end{lemma}

\begin{proof}
It suffices to prove (\ref{br-e11}), since (\ref{br-e10}) follows from it via
H\"older's inequality.

Let $p,p'\in(1,\infty)$ be dual exponents.
By H\"older's inequality, we have
\begin{align*}
\mathbb E \left[\Bigl(\int |\widehat{\mu}(\xi)|^{2} |\widehat{\mu}(-2\xi)| d\xi
\Bigr)^2\right] 
&\leq \mathbb E \bigl( ||\widehat{\mu}||_{2p}^4 ||\widehat{\mu}||_{p'}^2 \bigr)\\
&\leq \mathbb E \bigl( || \widehat{\mu}||_{2p}^{4p} \bigr) +
 \mathbb E \bigl( || \widehat{\mu}||_{p'}^{2p'} \bigr). \\
\end{align*}
Let $\eta,\kappa$ be positive numbers to be fixed later.  Using H\"older's
inequality once more, we estimate the last line by
\begin{align*}
&\mathbb{E}\left(\int|\whm(\xi)|^{4p}(1+|\xi|)^{1+2\eta}d\xi\right)
\left(\int(1+|\xi|)^{-(1+2\eta)}d\xi\right)\\
&\ \ +\mathbb{E}\left(\int|\whm(\xi)|^{2p'}(1+|\xi|)^{1+2\kappa}d\xi\right)
\left(\int(1+|\xi|)^{-(1+2\kappa)}d\xi\right)\\
&\leq C\int \mathbb{E}(|\whm(\xi)|^{4p})(1+|\xi|)^{1+2\eta}d\xi\\
&\ \ +C \int \mathbb{E}(|\whm(\xi)|^{2p'})(1+|\xi|)^{1+2\kappa}d\xi,
\end{align*}
where the constants are uniformly bounded as long as 
\begin{equation}\label{br-e20}
\eta,\kappa>\epsilon_0>0.
\end{equation}
By Proposition \ref{kahane-result}, this is bounded by 
\begin{align*}
& C\int (2C'_0p(1+|\xi|)^{-\alpha})^{2p}(1+|\xi|)^{1+2\eta}d\xi\\
& \ \ +C\int (C'_0p'(1+|\xi|)^{-\alpha})^{p'}(1+|\xi|)^{1+2\kappa}d\xi\\
& \leq C'\int (1+|\xi|)^{-2\alpha p+1+2\eta}d\xi
+ C'\int (1+|\xi|)^{-\alpha p'+1+2\kappa}d\xi.
\end{align*}
Both integrals above converge, provided that $\eta<\alpha p-1$ and $\kappa<
\frac{\alpha p'}{2}-1$. If we choose $p=3/2$, $p'=3$, then the conditions on
$\eta$ and $\kappa$ become $\eta,\kappa<\frac{3\alpha}{2}-1$, which is consistent
with (\ref{br-e20}) if $\alpha$ is bounded from below away from $2/3$.
\end{proof}
\begin{lemma}\label{br-lemma2}
There is an absolute constant $c_0>0$ such that
$\mathbb E(\Lambda(\mu, \mu, \mu)) > c_0$. 
\end{lemma}
\begin{proof}
Let us define 
\[ \Lambda_{\epsilon}(\mu, \mu, \mu) = \int \widehat{\mu}(\xi)^2
\widehat{\mu}(-2 \xi) e^{- 2\pi^2 \epsilon |\xi|^2} \, d\xi.\]
Lemma \ref{exp-conv} combined with dominated convergence (in $\xi$) yields $\Lambda_{\epsilon} \rightarrow \Lambda$ as $\epsilon \rightarrow 0$ for almost every Brownian path. Further, 
\begin{align*}
\mathbb E(\Lambda_{\epsilon})
&= \int \left[ \int \mathbb E \Bigl( e^{2 \pi i (2 W(t_3) - W(t_1) - W(t_2))\xi
- 2\pi^2 \epsilon |\xi|^2 }  \Bigr) d\xi \right] \prod_{i=1}^{3} d\theta(t_i)
\\ &= \sum_{\pi^{\ast}}  \int_{U_{\pi^{\ast}}} \left[ \int \mathbb E \Bigl(
e^{2 \pi i (2 W(t_3) - W(t_1) - W(t_2))\xi - 2\pi^2 \epsilon |\xi|^2 }  \Bigr)
d\xi \right] \prod_{i=1}^{3} d\theta(t_i),
\end{align*}
where $\pi^{\ast} = (\pi_1^{\ast}, \pi_2^{\ast}, \pi_3^{\ast})$ ranges over all permutations of $(1,2,3)$, and $U_{\pi^{\ast}} = \{t = (t_1, t_2, t_3) : t_{\pi_1^{\ast}} < t_{\pi_2^{\ast}} < t_{\pi_3^{\ast}} \}$. We claim that 
the sum above is strictly positive. 

Consider the term with $\pi^{\ast} = (1, 2, 3)$, all other cases being similar. Since $W(t_1)$, $W(t_2) - W(t_1)$ and $W(t_3) - W(t_2)$ are independent and normally distributed, $2 W(t_3) - W(t_1) - W(t_2) = (W(t_2) - W(t_1)) + 2(W(t_3) - W(t_2))$ is also normal with mean 0 and variance $\sigma^2 = \sigma_{\pi^{\ast}}^2(t) = t_2-t_1 + 4(t_3 - t_2)$. It follows therefore that 
\[ \mathbb E\bigl( e^{2 \pi i (2 W(t_3) - W(t_1) - W(t_2))\xi}\bigr)
= e^{-2 \pi^2 \sigma^2 |\xi|^2}, \text{ which implies}   \]
\begin{multline*} \int_{U_{\pi^{\ast}}} \left[ \int \mathbb E \Bigl(
e^{2 \pi i (2 W(t_3) - W(t_1) - W(t_2))\xi - 2\pi^2\epsilon |\xi|^2 }
 \Bigr) d\xi \right] \prod_{i=1}^{3} d\theta(t_i) \\ =
\frac{1}{\sqrt{2\pi}} \int_{U_{\pi}^{\ast}} \frac{\prod_{i=1}^{3}d\theta(t_i)}{\sqrt{(t_2-t_1) + 4(t_3-t_2) + \epsilon}}. \end{multline*}
Note that since $0\leq t_i\leq 1$, the expression in the denominator is bounded from above
by $\sqrt{5+\epsilon}$.

A similar calculation can be performed for each of the remaining terms.
Since $\Lambda_{\epsilon} \leq   \int |\widehat{\mu}(\xi)|^2 |\widehat{\mu}(-2 \xi)| d \xi$, whose expectation has been shown to be finite in Lemma \ref{exp-conv}, another application of dominated convergence (this time on the space of Brownian paths) yields 
\[ \mathbb E(\Lambda) = \lim_{\epsilon \rightarrow 0} \mathbb E(\Lambda_{\epsilon})
= \frac{1}{\sqrt{2\pi}} \int \left[\sum_{\pi^{\ast}} \frac{1_{U_{\pi^{\ast}}}(t)}{\sigma_{\pi^{\ast}}(t)} \right]\prod_{i=1}^3 d\theta(t_i). \]
Since the function $\sum_{\pi^{\ast}} 1_{U_{\pi^{\ast}}} \sigma_{\pi^{\ast}}^{-1}$
bounded from below by a strictly positive universal constant in $[0,1]^3$,
except for the zero-measure set where it is not defined, this proves the lemma.
\end{proof}

{\em Proof of Proposition \ref{prop-kahane2}.} The proposition follows
immediately from Lemmas \ref{exp-conv}, \ref{br-lemma2}, and the Paley-Zygmund
inequality \cite[p.8]{kahane}:
$$
\pp(X>\lambda \ee(X))\geq (1-\lambda)^2\frac{(\ee(X))^2}{\ee(X^2)},
$$
where $0<\lambda<1$ and $X$ is a positive random variable with $\ee(X^2)<\infty$.


\bibliographystyle{amsplain}

\noindent{\sc Department of Mathematics, University of British Columbia, Vancouver,
B.C. V6T 1Z2, Canada}
                                                                                     
\noindent{\it ilaba@math.ubc.ca, malabika@math.ubc.ca}

\end{document}